\documentclass[12pt, reqno, twoside]{amsart}
\usepackage{amsmath, amsthm, mathrsfs, amscd, amsfonts, amssymb, graphicx, color}
\usepackage[bookmarksnumbered, colorlinks, plainpages]{hyperref}
\hypersetup{colorlinks=true,linkcolor=red, anchorcolor=green, citecolor=cyan, urlcolor=red, filecolor=magenta, pdftoolbar=true}
\usepackage{enumitem}
\usepackage{enumitem}
\newtheorem{theorem}{Theorem}[section]
\newtheorem{lemma}[theorem]{Lemma}

\newtheorem{corollary}[theorem]{Corollary}
\theoremstyle{definition}

\theoremstyle{remark}
\newtheorem{remark}[theorem]{Remark}
\numberwithin{equation}{section}

\begin{document}

\setcounter{page}{1}

\title[Gradient estimates for nonlinear slow diffusion equations]
{Li-Yau Estimates and Harnack Inequalities for Nonlinear Slow Diffusion Equations on a Smooth Metric Measure Space}

\author[A. Taheri]{Ali Taheri}

\author[V. Vahidifar]{Vahideh Vahidifar}

\address{School of Mathematical and  Physical Sciences, 
University of Sussex, Falmer, Brighton, United Kingdom.}
\email{\textcolor[rgb]{0.00,0.00,0.84}{a.taheri@sussex.ac.uk}} 
 
\address{School of Mathematical and  Physical Sciences, 
University of Sussex, Falmer, Brighton, United Kingdom.}
\email{\textcolor[rgb]{0.00,0.00,0.84}{v.vahidifar@sussex.ac.uk}}

\subjclass[2020]{53C44, 58J60, 58J35, 60J60}

\keywords{Smooth metric measure space, Nonlinear slow diffusion equation, 
Aronson-B\'enilan estimates, Differential Harnack inequalities, Li-Yau gradient estimates, Harnack inequalities}


\begin{abstract}
We present new gradient estimates and Harnack inequalities for positive solutions to nonlinear slow diffusion 
equations. The framework is that of a smooth metric measure space 
$(\mathscr M,g,d\mu)$ with invariant weighted measure $d\mu=e^{-\phi} dv_g$ and diffusion operator 
$\Delta_\phi=e^\phi {\rm div} (e^{-\phi} \nabla)$ -- the $\phi$-Laplacian. The nonlinear slow diffusion 
equation, then, for $x \in {\mathscr M}$ and $t>0$, and fixed exponent $p>1$, takes the form 
\begin{equation*}
\partial_t u (x,t) - \Delta_\phi u^p (x,t) = \mathscr N (t,x,u(x,t)).
\end{equation*}
We assume that the metric tensor $g$ and potential $\phi$ are space-time dependent; hence the same is true of the 
usual metric and potential dependent differential operators and curvature tensors. 
The estimates are established under natural lower bounds on the Bakry-\'Emery $m$-Ricci curvature tensor and the time 
derivative of metric tensor. The curious interplay between geometry, nonlinearity and evolution and their influence 
on the estimates is at the centre of this investigation. The results here considerably extend and improve earlier results 
on slow diffusion equations. Several implication, special cases and corollaries are presented and discussed. 
\end{abstract}

\maketitle
 
{ 
\hypersetup{linkcolor=black}
\tableofcontents
}

\allowdisplaybreaks

\section{Introduction} 
\label{sec1}

In this paper we are concerned with gradient estimates of Li-Yau or Aronson-B\'enilan type (also known as differential 
Harnack estimates) for positive solutions to nonlinear slow diffusion equations. Here, at the start, we mention that 
throughout the paper, the terms slow diffusion equation and porous medium equation will be used synonymously. 
Such estimates are known to be of huge utility and significance in geometric analysis. For example, they play a profound 
role in our understanding of the short and long term dynamics, the interaction between geometry and curvature on the 
one hand and diffusion on the other, and the derivation of various qualitative and structural properties of solutions 
such as parabolic Harnack inequalities, Liouville-type results, ancient and eternal solutions, 
blowup profiles and heat kernel bounds to mention a few.

The framework we study in this paper is that of a smooth metric measure space or a weighted 
manifold. This is a triple $(\mathscr M,g,d\mu)$ where $({\mathscr M}, g)$ is a complete 
Riemannian manifold of dimension $n \ge 2$, $d\mu=e^{-\phi} dv_g$ is a weighted measure associated with 
the metric $g$ and potential $\phi$ and $dv_g$ is the usual Riemannain volume measure on ${\mathscr M}$. 
The nonlinear slow diffusion equation for a space-time function $u=u(x,t)$ can then be written as
\begin{align} \label{eq11}
\partial_t u(x,t) - \Delta_\phi u^p(x,t) = \mathscr N(t,x,u(x,t)). 
\end{align}
Here $p>1$ is a fixed exponent, $\Delta_\phi$ is the Laplace operator associated with $({\mathscr M}, g, d\mu)$ 
(called the Witten or weighted or $\phi$-Laplacian) defined for $w \in \mathscr{C}^2(\mathscr M)$ by 
\begin{equation} \label{f-Lap-definition}
\Delta_\phi w = {\rm div}_\phi (\nabla w) = e^\phi {\rm div}(e^{-\phi} \nabla w) = \Delta w - \langle \nabla \phi, \nabla w\rangle,
\end{equation} 
and ${\mathscr N}={\mathscr N}(t,x,u)$ 
is a sufficiently smooth nonlinearity (or forcing term) depending on time $t$, the spatial variable $x$,  
and the dependent variable $u$. We refer the reader to \cite{Aron, DaskKe, Vaz} and 
the references therein for background on slow diffusion equations (i.e., porous medium equations), 
the underlying theory and the numerous applications.

The study of gradient estimates for diffusion equations on Riemannian manifolds goes back to the seminal paper \cite{[LY86]}. 
Subsequent advances and developments, in particular, to the context of evolving manifolds and/or smooth metric measure spaces came later 
({\it see} \cite{AM, BaCP, Chow, Ha93, SZ, Taheri-GE-1, Taheri-GE-2, TVahNA, TVahCurv, TV-FDE-HS-a, TVDiffHar, Wu18, WuW14, Zhang}). 
The study of gradient estimates for degenerate and nonlinear diffusion equations is more recent. For the slow diffusion 
equation on Riemannian manifolds we can refer to the works \cite{AronBen,Barenblatt,BBGM, BDM, Bon, Huang-Huang-Li, Huang, PLu, MZS, Zhu} and for the 
slow diffusion equation on smooth metric measure spaces we can refer to the works \cite{ Hunag-Li, TV-PME-a}. 
\footnote{We discuss the fast diffusion equation corresponding to $0<p<1$ in \eqref{eq11} in a forthcoming paper.}

The principal aim of this paper is to push the subject forward by 
establishing gradient estimates for the nonlinear slow diffusion equation \eqref{eq11} in the context of evolving triples 
$({\mathscr M},g,d\mu)$. Such problems have particularly moved to the forefront of research in geometric analysis since 
the work of G.~Perelman \cite{Pe02} and the study of forward and backward heat-type equations and gradient flows on 
triples $({\mathscr M},g,d\mu)$ evolving under (super) Hamilton-Ricci or (super) 
Perelman-Ricci flows ({\it see} \cite{AGS, BaCP, Chow, MReto, Sturm, Taheri-GE-1, Taheri-GE-2, TVahNA, TVDiffHar, VC, Zhang}). 
As can be readily seen the addition of evolution and nonlinearity brings in considerable technical challenges to 
the landscape but naturally has huge applications and implications.

Let us elaborate more on the previous two paragraphs to put the problem and our results in better perspective. 
Starting with an $n$-dimensional complete Riemannian manifold $({\mathscr M}, g)$, Li and Yau in \cite{[LY86]} 
studied positive solutions to the heat equation 
\begin{equation}\label{EQ-1.1}
\partial_t u = \Delta u. 
\end{equation}
\textbf{\bf Theorem A.} (Li-Yau \cite{[LY86]}) Let $({\mathscr M}, g)$ be a complete Riemannian manifold with 
${\mathscr Ric}(g) \ge -{\mathsf k}g$ in ${\mathscr B}_{2R} \subset {\mathscr M}$ for some ${\mathsf k} \ge 0$. 
Suppose that $u$ is a positive solution to \eqref{EQ-1.1} on ${\mathscr B}_{2R} \times [0,T]$. Then for every 
$x$ in ${\mathscr B}_R$ and $0<t \le T$ we have
\begin{equation}\label{EQ-1.2}
\frac{|\nabla u|^2}{u^2} - \alpha \frac{u_t}{u} \le \frac{C\alpha^2}{R^2} 
\left[ \frac{\alpha^2}{\alpha - 1} + \sqrt{{\mathsf k} R} \right] + \frac{n\alpha^2 {\mathsf k}}{2(\alpha - 1)} + \frac{n\alpha^2}{2t},
\end{equation}
where $\alpha > 1$ is a constant and $C>0$ depends only on $n$. Letting $R \nearrow \infty$, 
\eqref{EQ-1.2} gives the following estimate on the complete noncompact Riemannian manifold $({\mathscr M}, g)$:
\begin{equation}\label{EQ-1.3}
\frac{|\nabla u|^2}{u^2} - \alpha \frac{u_t}{u} \le \frac{n\alpha^2 {\mathsf k}}{2(\alpha - 1)} + \frac{n\alpha^2}{2t}.
\end{equation}

Later on E.B.~Davies in \cite{Davies} improved the global estimate \eqref{EQ-1.3} under the same conditions to the following:
\qquad \\ 
\textbf{Theorem B.} (Davies \cite{Davies}) Let $({\mathscr M}, g)$ be a complete noncompact Riemannian manifold 
with ${\mathscr Ric}(g) \ge -{\mathsf k}g$ for some ${\mathsf k} \ge 0$. Suppose that $u$ is a positive solution to \eqref{EQ-1.1}. Then
\begin{equation}\label{EQ-1.4}
\frac{|\nabla u|^2}{u^2} - \alpha \frac{u_t}{u} \leq \frac{n\alpha^2 {\mathsf k}}{4(\alpha - 1)} + \frac{n\alpha^2}{2t}.
\end{equation}

Interestingly, by using the idea of matrix Harnack inequalities, a little later R.~Hamilton in \cite{Ha93} proved the following:
\qquad \\
\textbf{Theorem C.} (Hamilton \cite{Ha93}) Let $({\mathscr M}, g)$ be a complete noncompact Riemannian manifold with 
${\mathscr Ric}(g) \ge -{\mathsf k}g$ for some $k \ge 0$. Suppose that $u$ is a positive solution to \eqref{EQ-1.1}. Then
\begin{equation}
\frac{|\nabla u|^2}{u^2} - e^{2{\mathsf k}t} \frac{u_t}{u} \leq e^{4{\mathsf k}t} \frac{n}{2t}. \tag{1.5}
\end{equation}

Whilst the above estimates and improvements relate to the heat equation \eqref{EQ-1.1}, the simplest nonlinear generalisation 
of  \eqref{EQ-1.1} and the equation that is closer to the subject matter of this paper is the slow diffusion equation ($p>1$)
\begin{equation}\label{EQ-1.6}
\partial_t u = \Delta u^p.
\end{equation}
\textbf{Theorem D.} (Lu, Ni, Vázquez, Villani \cite{PLu}) Let $({\mathscr M}, g)$ be a complete Riemannian 
manifold with ${\mathscr Ric}(g) \ge -{\mathsf k}g$ in ${\mathscr B}_{2R}\subset {\mathscr M}$ for some ${\mathsf k} \ge 0$. 
Suppose that $u$ is a positive solution to \eqref{EQ-1.6}. Let $v = p u^{p-1}/ (p-1)$ and $M = \sup v$. 
Then for any $\alpha > 1$, on the ball $B_p(R)$, we have
\begin{equation}\label{EQ-1.7}
\frac{|\nabla v|^2}{v} - \alpha \frac{v_t}{v} 
\le b \alpha^2 \left\{ \frac{CM}{R^2} \left[ \frac{p^2 b\alpha^2}{\alpha-1} + (p-1)(1 + \sqrt{{\mathsf k}} R) \right] 
+ \frac{p-1}{\alpha-1} M {\mathsf k} + \frac{1}{t} \right\},
\end{equation}
where $b = [n(p-1)]/[n(p-1)+2]$ and the constant $C$ depends only on $n$. Letting 
$R \nearrow \infty$, \eqref{EQ-1.7} gives the following estimate on the complete 
noncompact Riemannian manifold $({\mathscr M}, g)$:
\begin{equation}
\frac{|\nabla v|^2}{v} - \alpha \frac{v_t}{v} \le \frac{b \alpha^2(p-1)}{\alpha-1} M {\mathsf k} + \frac{b \alpha^2}{t}.
\end{equation}
Some improvements and generalisations of this result can be found in \cite{Huang-Huang-Li, Hunag-Li, Huang}. 
As is seen \eqref{EQ-1.6} is a simple and special case of the general nonlinear porous medium equation \eqref{eq11} 
that is treated in this paper. Our aim here is to establish gradient estimates for \eqref{eq11} and investigate closely how evolution, 
nonlinearity and the geometry of $({\mathscr M},g,d\mu)$ influence the estimate and hence 
its consequences, in particular, the parabolic Harnack inequalities and Liouville results.

\qquad \\
{\bf Plan of the paper.} The main results are presented in Section \ref{sec2}. Here we formulate the two 
main estimates in their local and global forms as well as the resulting Harnack inequalities for positive 
solutions to \eqref{eq11}. In doing so we introduce various quantities pertaining to the geometry and 
evolution of the triple $({\mathscr M}, g, d\mu)$, its bounds and the nonlinearity ${\mathscr N}={\mathscr N}(t,x,u)$. 
We discuss several implication and consequences of the estimates, in particular, to the static triple case. The rest 
of the paper is devoted to the proof of the main results. In Section \ref{sec3} we introduce a Harnack quantity 
depending on the solution $u$ and describe its evolution.  In Section \ref{sec4} we use the results in Section \ref{sec3} 
to prove the local estimate in Theorem \ref{thm5.1}. In Section \ref{sec5} we present the proof of the Harnack inequality 
in Corollary \ref{paraHar}. In Section \ref{sec6} we present the proof of the  local estimate in Theorem \ref{Thrm-2.5}. 
Here again we use the results in Section \ref{sec3}. In Sections \ref{sec4} and \ref{sec6} localisation and 
construction of suitable cut-off functions plays a key role. Finally, in Section \ref{sec7} we present the proof of 
the Harnack inequality in Corollary \ref{corollary-Harnack-global.n}.

\qquad \\
{\bf Notation.} 
Fixing a reference point $x_0 \in {\mathscr M}$ we denote by $d=d(x,x_0, t)$ the Riemannian distance between 
$x$ and $x_0$ at time $t$ on ${\mathscr M}$ with respect to the metric $g=g(t)$. We write $\varrho=\varrho(x,x_0,t)$ 
for the geodesic radial variable measuring distance between $x$ and origin $x_0$ at time $t>0$. For $R>0$, 
$T>0$ we define the compact space-time cylinder 
\begin{equation}
Q_{R,T} (x_0)\equiv \{ (x, t) | d(x, x_0, t) \le R, 0 \le t \le T \} \subset {\mathscr M} \times [0, T].
\end{equation}
When $g$ is time independent, we denote by $\mathscr{B}_R(x_0) \subset {\mathscr M}$ the 
geodesic ball of radius $R>0$ centred at $x_0$. It is evident that in this case we have 
\begin{equation}
Q_{R,T} (x_0) = \mathscr{B}_R(x_0) \times [0, T] \subset {\mathscr M} \times [0, T].
\end{equation}
When the choice of $x_0$ is clear from the context we often abbreviate and write $d(x, t)$, 
$\varrho(x,t)$ or ${\mathscr B}_R$, $Q_{R,T}$ respectively.

For a function of multiple variables, we denote its partial derivatives with subscripts accordingly. 
For instance for  $\Gamma=\Gamma(x,v)$ we denote its partial derivatives with respect to 
$x=(x_1, \dots, x_n)$ or $v$ by  $\Gamma_x$ or $\Gamma_v$ respectively. Moreover we reserve 
the notation $\Gamma^x$ for the function $x \mapsto \Gamma (x, v)$ obtained by freezing the 
argument $v$ and viewing it as a function of $x$. In the sequel we frequently make use of the 
notations $\nabla \Gamma^x$ and $\Delta_f \Gamma^x$.  Finally we write $z=z_+ + z_-$ with 
$z_+=\max(z, 0)$ and $z_-=\min(z, 0)$.

The Bakry-\`Emery $m$-Ricci curvature tensor ${\mathscr Ric}^m_\phi(g)$ associated with $({\mathscr M}, g, d\mu)$ 
is defined as
\begin{equation} \label{Ricci-m-f-intro}
{\mathscr Ric}^m_\phi(g) := {\mathscr Ric}(g) + \nabla\nabla \phi - \frac{\nabla \phi \otimes \nabla \phi}{m-n}, \qquad m \ge n. 
\end{equation}
Here ${\mathscr Ric}(g)$ is the usual Riemannain Ricci curvature tensor of $g$, $\nabla \nabla \phi={\rm Hess}(\phi)$ 
denotes the Hessian of $\phi$, and $m \ge n$ is a constant ({\it see} \cite{Bak}). For the sake of clarification, in 
\eqref{Ricci-m-f-intro}, when $m=n$, by convention $\phi$ is only allowed to be a constant, thus giving 
${\mathscr Ric}^m_\phi(g)={\mathscr Ric}(g)$. We also recall the weighted Bochner-Weitzenb\"ock formula asserting 
that for any function $w$ of class ${\mathscr C}^3(M)$ we have, 
\begin{equation} \label{Bochner-1}
\frac{1}{2} \Delta_\phi |\nabla w|^2 - \langle \nabla w, \nabla \Delta_\phi w \rangle
= |\nabla\nabla w|^2 + {\mathscr Ric}_\phi (\nabla w, \nabla w), 
\end{equation}
where ${\mathscr Ric}_\phi(g)={\mathscr Ric}(g) + \nabla\nabla \phi$, that is, formally, 
${\mathscr Ric}_\phi(g) = \lim {\mathscr Ric}^m_\phi(g)$ as $m \nearrow \infty$.
The Riemannian case is the special case $\Delta_\phi=\Delta$, 
${\mathscr Ric}_\phi^m(g) \equiv {\mathscr Ric}(g)$ and $m=n$.


\section{Statement of the main results}
\label{sec2}

In this section we present the main results. The proofs will appear in later sections 
after developing the necessary tools. Note that since the metric $g$ and potential $\phi$ here 
are time dependent, then the same is true of the weighted measure $d\mu=e^{-\phi} dv_g$, the 
differential operators ${\rm div}$, $\nabla$, $\Delta$, $\Delta_\phi$ and the 
$m$-Ricci curvature tensor 
\begin{equation}
{\mathscr Ric}_\phi^m(g)(x,t) = {\mathscr Ric} (g)(x,t) + \nabla \nabla \phi (x,t) 
- \dfrac{\nabla \phi(x,t) \otimes \nabla \phi(x,t)}{m-n}. 
\end{equation}

In order to proceed with the estimates we need to formulate suitable bounds on the geometry and evolution 
of the triple $({\mathscr M}, g, d\mu)$. To this end we start by setting  
\begin {align} \label{2.4}
\frac{\partial g}{\partial t} (x,t)= 2h (x,t), \qquad (x,t) \in {\mathscr M} \times [0, T]. 
\end{align}
With the aid of this notation we can pose certain pointwise bounds on the metric tensor $g$ and 
potential $\phi$ in suitable space-time regions as needed later. These take the forms:  
\begin{align} \label{EQ-eq-2.2}
{\mathscr Ric}_\phi^m (g) (x,t) \ge -(m-1)k g (x,t),
\end{align}
\begin{align} \label{vahideh-A-2} 
-\underline k_1 g(x,t) \le h(x,t) \le \overline k_1 g(x,t), 
\end{align}
\begin{align} \label{vahideh-A-2-b}
|\nabla h (x,t)| \le k_2,
\end{align}
\begin{align} \label{vahideh-A-3}
|\nabla \phi (x,t)| \le \ell_1, \qquad |\nabla \partial_t \phi(x,t)| \le \ell_2.  
\end{align}
Here $k, \underline k_1, \overline k_1, k_2 \ge 0$ and $\ell_1, \ell_2 \ge 0$ are appropriate constants. 
We note that \eqref{EQ-eq-2.2}--\eqref{vahideh-A-2} are to be understood in the sense 
of symmetric second order tensors. Throughout the paper we write $b = [m(p-1)]/[1+ m(p-1)]$ 
(see Remark \ref{remark-on-b} for more on this). In view of $p>1$ and $m \ge n \ge 2$ we 
have $0<b<1$. Note that $m$ is a real number (not necessarily an integer). 
Unless otherwise stated we  assume $\alpha= \alpha(t) > 1$ nondecreasing and $\beta= \beta(t)$ 
are   functions of class $\mathscr{C}^1[0,T]$ with $T>0$ fixed as in \eqref{2.4}. 
See Remark \ref{interesting-alpha-beta-remark} for some special examples of interest. 
Throughout the paper we make use of the pressure transform $v =pu^{p-1}/(p-1)$, relating 
$u>0$ to a corresponding $v>0$. In this respect it is convenient to introduce the rescaled 
nonlinearity $\mathscr G= \mathscr G(t,x,v)$ by writing 
\begin{align}\label{eq5.2}
\mathscr G(t,x,v)= p [(p-1)v/p]^{\frac{p -2}{p-1}} 
\mathscr N \left(t,x,[(p-1)v/p]^{\frac{1}{p-1}}\right). 
\end{align}

\qquad \\
{\bf The first main gradient estimate.} 
We start by introducing a string of non-negative quantities depending on the geometry and evolution of the triple 
$({\mathscr M}, g, d\mu)$. Indeed, using the constants in \eqref{EQ-eq-2.2}--\eqref{vahideh-A-3} 
we introduce for $p>1$, $c_1, R>0$ and $M\ge0$, 
\begin{align}
\mathsf K & = 2c_1 \underline k_1 + c_1^2 b \alpha^2 p^2 M /[2(\alpha -1)R^2], \label{EQ-eq-K}\\ 
\mathsf L &= \alpha (p-1)  \ell_2/2 + \alpha (p-1) \underline k_1 \ell_1, \label{EQ-eq-L}\\
\mathsf M &= \alpha^2 (p-1)n[(\underline k_1 + \overline k_1)^2+2k_2], \label{EQ-eq-M}\\
\mathsf N &= 2(p-1)M[(m-1) k+k_2] +2 (\alpha -1) \overline k_1, \label{EQ-eq-N}
\end{align}
and for $0<\varepsilon<2(\alpha-1)^2/(b \alpha^2)$ on $[0, T]$, 
\begin{align}
\mathsf E = (3/2)^{3/2} M/ \sqrt \varepsilon, \qquad 
\mathsf F = b \alpha^2 /[4 (\alpha -1)^2-2 \varepsilon b \alpha^2]. \label{EQ-eq-E}
\end{align}
We note that these quantities in general depend on $t$ through $\alpha$, unless, of course, the latter 
is constant. Associated with ${\mathscr N}$ in \eqref{eq11}, $v=v(x,t)>0$ and 
$Q_{R, T} \subset {\mathscr M} \times [0, T]$, $M=M(R)=\sup v$ over $Q_{R,T}$ and with the aid of \eqref{EQ-eq-K}--\eqref{EQ-eq-E} 
we further introduce 
\begin{align}
&\mu_1(R) =\sup_{Q_{R,T}} \left [\mathscr G_v +\frac{\alpha'}{\alpha} -\frac{2 \beta}{b \alpha^2} +\mathsf K\right]_+,\label{V-A-1}\\
&\mu_2(R) = \sup_{Q_{R,T}}\left[ \sqrt {\mathsf E} \left\{ (\alpha-1) \frac{|\mathscr G_x|}{v}+ \alpha (p-1)  |\mathscr G_{x v}|+ \mathsf L \right\} \right],\label{V-A-2}\\
&\mu_3(R) = \sup_{Q_{R,T}} \left[ \beta \mathscr G_v - \alpha (p-1) \Delta_\phi \mathscr G^x 
+\frac{\beta}{\alpha}\left( \alpha' - \frac{\beta}{b \alpha} \right) - \beta' + \mathsf M\right]_+,\label{V-A-3}\\
&\mu_4(R) = \sup_{Q_{R,T}} \left[ \sqrt \mathsf F \left\{ (\alpha -1)\left(\frac{\mathscr G}{v} -\mathscr G_v \right)
- \alpha(p-1) v\mathscr G_{vv} - \frac{2(\alpha -1)}{b \alpha^2} \beta
- \frac{\alpha'}{\alpha}+\mathsf N \right\}\right] _+. \label{V-A-4}
\end{align}
With this in place we are now in a position to state the first main result of the paper.

\begin{theorem}\label{thm5.1}
Let $(\mathscr M, g, d\mu)$ be a complete smooth metric measure space with $d\mu=e^{-\phi} dv_g$. 
Suppose that the metric and potential are time dependent, of class $\mathscr{C}^2$ and that for suitable 
constants $k, \underline k_1, \overline k_1, k_2\ge 0$, $\ell_1, \ell_2 \ge 0$ and $n \le m < \infty$ satisfy 
the bounds \eqref{EQ-eq-2.2}--\eqref{vahideh-A-3} in the space-time cylinder $Q_{2R,T}$. 
Let $u$ be a positive solution to equation \eqref{eq11} with $p>1$ and set $v =pu^{p-1}/(p-1)$. Then for 
any $0<\varepsilon<2[\alpha(t)-1]^2/[b \alpha^2(t)]$ and every $(x,t)$ in $Q_{R,T}$ 
with $t>0$ we have 
\begin{align} \label{estimate-5.4}
\left[\frac{|\nabla v|^2}{\alpha v} -\frac{\partial_t v}{v} + \frac{\mathscr G}{v} -\frac{\beta}{\alpha} \right](x,t)
\le&~ \frac{b \alpha}{t} + b \alpha \mu_1(2R) \\
&+ \sqrt b \bigg\{\mu_2^{4/3}(2R) + \mu_3(2R) + \mu_4^2(2R)\bigg\}^{1/2}\nonumber\\
&+b (p-1)\alpha \frac{M(2R)}{R^2}[c_2 +(m-1) c_1(1+R \sqrt{k})+ 2c_1^2].\nonumber
\end{align}
Here $M = \sup v$ over $Q_{2R,T}$, $c_1$, $c_2 >0$ are the constants in \eqref{9.30} 
and the non-negative quantities $\mu_j$ $($$1 \le j \le 4$$)$ are as in \eqref{V-A-1}--\eqref{V-A-4}. 
\end{theorem}

Subject to the bounds in Theorem \ref{thm5.1} being global (in space) by passing to the limit 
$R \nearrow \infty$ we obtain the following global (in space) version of the estimate in 
Theorem \ref{thm5.1}.

\begin{corollary}\label{corollary-A-Harnack}
If the bounds \eqref{EQ-eq-2.2}--\eqref{vahideh-A-3} hold globally on $\mathscr M \times [0, T]$
and $v$ is bounded on $\mathscr M \times [0, T]$ 
then for every $x \in \mathscr M$ and $0< t<T$ we have 
\begin{align}\label{estimate-5.4-global}
\left[\frac{|\nabla v|^2}{\alpha v} -\frac{\partial_t v}{v} + \frac{\mathscr G}{v} -\frac{\beta}{\alpha} \right](x,t)
\le&~ \frac{b \alpha}{t} + b \alpha \mu_1 
+ \sqrt b \bigg\{\mu_2^{4/3} + \mu_3 + \mu_4^2\bigg\}^{1/2}.
\end{align}
Here $M = \sup v$ over $\mathscr M \times [0, T]$
and $\mu_j \ge 0$ $($$1 \le j \le 4$$)$ are as in \eqref{V-A-1}--\eqref{V-A-4}
with $\mathsf K= 2c_1 \underline k_1$ in $\mu_1$ and all the supremums being taken over $\mathscr M \times [0, T]$.
\end{corollary}

One of the immediate applications of the above estimates is to parabolic Harnack inequalities. 
Here we present the formulation in the global case with $\alpha>1$ constant. 
The local version of the inequality is analogous and will be abbreviated.

\begin{corollary} \label{paraHar}
Under the assumptions of Corollary $\ref{corollary-A-Harnack}$, for every constant $\alpha >1$, 
and for all $x_1, x_2 \in \mathscr M$ and $0 <t_1<t_2<T$ we have 
\begin{align}\label{Eq5.62Har}
\frac{v(x_1, t_1)}{v(x_2, t_2)} \le 
{\rm exp} \left[ \frac{\alpha \mathsf L(x_1,x_2)}{4 \underline M (t_2-t_1)}
+ \frac{\mathsf H}{\alpha} (t_2 -t_1)\right] \left(\frac{t_2}{t_1}\right)^{b \alpha}.
\end{align}
Here $\underline M = \inf v$ over $\mathscr M \times [0, T]$, $\mathsf L=\mathsf L(x_1,x_2)$ is the non-negative quantity 
\begin{align} \label{L-Harnack-global-equation}
\mathsf L(x_1,x_2) = \inf_{\gamma \in \Gamma(x_1, x_2)} ||\,|\dot \gamma(t)|_{g(t)}^2 ||_{L^1(0,1; dt)},
\end{align} 
where $\Gamma (x_1, x_2) =\{ \gamma \in \mathscr{C}^1( [0, 1]; \mathscr M) : \gamma(0) = x_1, \gamma(1) = x_2\}$, 
and $\mathsf H$ is given by  
\begin{align}\label{V-A-H}
\mathsf H=&~ \mu_0 +b \alpha^2 \mu_1
+ \alpha \sqrt b \bigg\{\mu_2^{4/3} + \mu_3 + \mu_4^2\bigg\}^{1/2},
\end{align}
with $\mu_0 = \sup [\beta -\alpha \mathscr G/v]$ over ${\mathscr M} \times [0, T]$.
The constants in \eqref{V-A-H} are as in \eqref{estimate-5.4-global}.
\end{corollary}

Before moving on to presenting the second estimate we pause briefly to discuss some special cases 
and implications of the estimates in Theorem \ref{thm5.1} and Corollary \ref{corollary-A-Harnack}.   
\begin{itemize}
\item When $\mathscr N = \mathscr N(u)$ we have $\mathscr G (v) = p [(p-1)v/p]^{\frac{p -2}{p-1}} 
\mathscr N ([(p-1)v/p]^{\frac{1}{p-1}})$. In this case the terms ${\mathscr G}_x$, ${\mathscr G}_{xv}$ 
and $\Delta_\phi {\mathscr G}^x$ in $\mu_2(R)$ and $\mu_3(R)$ are zero. Moreover, if $\alpha >1$ is 
constant and $\beta \equiv 0$ we have $\mu_2(R) = \sqrt{\mathsf E} \mathsf L$, $\mu_3(R) = \mathsf M$ and 
$\mu_1(R), \mu_4(R)$ simplify to 
\begin{align}\label{EQ-2.16}
\mu_1(R) &= \sup_{Q_{R,T}} [\mathscr G_v+ \mathsf K]_+,\\
\mu_4(R) &= \sqrt {\mathsf F}\sup_{Q_{R,T}} [(\alpha -1) (\mathscr G/v -\mathscr G_v)
- \alpha(p-1) v\mathscr G_{vv} +\mathsf N]_+.
\end{align}
\item In the context described above if $\mathscr G_v \le 0$ and 
$\alpha(p-1) v\mathscr G_{vv}- (\alpha -1) (\mathscr G/v -\mathscr G_v) \ge 0$ 
then $\mu_1(R) \le \mathsf K$ and $\mu_4(R) \le \sqrt{\mathsf F} \mathsf N$. 
In this case the local estimate \eqref{estimate-5.4} gives
\begin{align}\label{EQ-eq-2.17}
\frac{|\nabla v|^2}{\alpha v} -\frac{\partial_t v}{v} + \frac{\mathscr G(v)}{v}
\le&~ \frac{b \alpha}{t} + b \alpha \mathsf K 
+ \sqrt b \bigg\{ (\sqrt {\mathsf E} {\mathsf L})^{4/3}
+ \mathsf M + \mathsf F \mathsf N^2 \bigg\}^{1/2}\\
&+b (p-1)\alpha \frac{M(2R)}{R^2}[c_2 +(m-1) c_1(1+R \sqrt{k})+ 2c_1^2 ], \nonumber 
\end{align}
and the global estimate \eqref{estimate-5.4-global} gives 
\begin{align}\label{EQ-eq-2.18}
\frac{|\nabla v|^2}{\alpha v} -\frac{\partial_t v}{v} + \frac{\mathscr G(v)}{v}
\le&~b \alpha \left\{\frac{1}{t} + 2c_1 \underline k_1\right\} \nonumber \\
&+ \sqrt b \bigg\{ (\sqrt \mathsf E \mathsf L)^{4/3}
+ \mathsf M + \mathsf F \mathsf N^2 \bigg\}^{1/2}.
\end{align}
\item When $\mathscr N \equiv 0$ (equivalently $\mathscr G \equiv 0$), $\alpha >1$ 
is a constant and $\beta \equiv 0$ we have 
\begin{align}
\mu_1(R) = \mathsf K,\quad \mu_2(R)=\sqrt{\mathsf E} \mathsf L, \quad
\mu_3(R)=\mathsf M, \quad \mu_4(R)= \sqrt{\mathsf F} \mathsf N. 
\end{align}
In this case the local estimate \eqref{estimate-5.4} reduces to \eqref{EQ-eq-2.17}
and the global estimate \eqref{estimate-5.4-global} reduces to \eqref{EQ-eq-2.18}. 
(In either case without the term $\mathscr G/v$ on the left.) 
\end{itemize}

\qquad \\
{\bf The second main gradient estimate.}
We now move on to the second main estimate in the paper. The formulation and notation 
is mostly similar to those used in Theorem \ref{thm5.1}, however, there are some major 
differences that we describe below. Indeed, with a slight variation of 
\eqref{EQ-eq-M}--\eqref{EQ-eq-E}, 
for $\alpha=\alpha(t)>1$ and $\beta=\beta(t)$ as before we set 
\begin{align}
\tilde {\mathsf M} &= (p-1)n[\alpha(\underline k_1 + \overline k_1)^2+2 k_2], \label{EQ-eq-M.n} \\
\tilde{\mathsf N} &= 2(p-1)M[(m-1) k/\alpha+ k_2] +2 (\alpha -1) \overline k_1/\alpha, \label{EQ-eq-N.n}\\
\tilde{\mathsf F} &= b \alpha^3 /[4 (\alpha -1)^2-2 \varepsilon b \alpha^3]. \label{EQ-eq-E.n}
\end{align}
Moreover, in analogy with \eqref{V-A-1}--\eqref{V-A-4} and by recalling \eqref{EQ-eq-K}--\eqref{EQ-eq-L}, 
we introduce the non-negative quantities 
\begin{align}
&\lambda_1(R) =\sup_{Q_{R,T}} \left [ \mathscr G_v -\frac{2 \beta}{b \alpha^2}+\mathsf K \right]_+,\label{V-A-1.n}\\
&\lambda_2(R) = \sup_{Q_{R,T}}\left[ \sqrt {\mathsf E} \left\{ \frac{\alpha-1}{\alpha} \frac{|\mathscr G_x|}{v}+ (p-1) |\mathscr G_{x v}|
+ \frac{\mathsf L}{\alpha} \right\} \right],\label{eq-6.29}\\
&\lambda_3(R) = \sup_{Q_{R,T}} \left[ \frac{\beta}{\alpha} \mathscr G_v - (p-1) \Delta_\phi \mathscr G^x 
+\frac{\beta}{\alpha^2}\left( \alpha' - \frac{\beta}{b \alpha} \right) - \frac{\beta'}{\alpha} 
+ \tilde {\mathsf M} \right]_+,\label{eq-6.30}\\
&\lambda_4(R) = \sup_{Q_{R,T}} \left[ \sqrt{\tilde{\mathsf F}} \left\{ \frac{\alpha -1}{\alpha} 
\left(\frac{\mathscr G}{v} -\mathscr G_v \right)- (p-1) v\mathscr G_{vv} 
- \frac{2(\alpha -1)}{b \alpha^3} \beta- \frac{\alpha'}{\alpha^2} + \tilde{\mathsf N} \right\} \right]_+. \label{V-A-4.n}
\end{align}

\begin{theorem}\label{Thrm-2.5}
Let $(\mathscr M, g, d\mu)$ be a complete smooth metric measure space with $d\mu=e^{-\phi} dv_g$. 
Suppose that the metric and potential are time dependent, of class $\mathscr{C}^2$ and that for suitable 
constants $k, \underline k_1, \overline k_1, k_2\ge 0$, $\ell_1, \ell_2 \ge 0$ and $n \le m < \infty$ satisfy 
the bounds \eqref{EQ-eq-2.2}--\eqref{vahideh-A-3} in the space-time cylinder $Q_{2R,T}$. 
Let $u$ be a positive solution to equation \eqref{eq11} with $p>1$ and set $v =pu^{p-1}/(p-1)$. Then for 
any $0<\varepsilon<2[\alpha(t)-1]^2/[b \alpha^3(t)]$ and every $(x,t)$ in $Q_{R,T}$ 
with $t>0$ we have
\begin{align} \label{estimate-2-local}
\left[\frac{|\nabla v|^2}{\alpha v} -\frac{\partial_t v}{v} 
+ \frac{\mathscr G}{v} -\frac{\beta}{\alpha} \right](x,t)
\le&~\frac{b \alpha}{t} + b \alpha \lambda_1(2R) \\
&+\sqrt {b \alpha} \bigg\{\lambda_2^{4/3}(2R) + \lambda_3(2R) + \lambda_4^2(2R)\bigg\}^{1/2} \nonumber\\
&+ b (p-1)\alpha \frac{M(2R)}{R^2} [c_2 +(m-1) c_1(1+R \sqrt{k})+ 2c_1^2]. \nonumber
\end{align}
Here $M = \sup v$ over $Q_{2R,T}$, $c_1$, $c_2 >0$ are the constants in \eqref{9.30} 
and the non-negative quantities $\lambda_j$ $($$1 \le j \le 4$$)$ are as in \eqref{V-A-1.n}--\eqref{V-A-4.n}. 
\end{theorem}

Again subject to the bounds in Theorem \ref{Thrm-2.5} being global (in space) by passing to the limit 
$R \nearrow \infty$ we obtain the following global (in space) version of the estimate.

\begin{corollary} \label{corollary-global.n}
If the bounds \eqref{EQ-eq-2.2}--\eqref{vahideh-A-3} hold globally on $\mathscr M \times [0, T]$
and $v$ is bounded on $\mathscr M \times [0, T]$ 
then for every $x \in \mathscr M$ and $0< t<T$ we have 
\begin{align}\label{estimate-2-global}
\left[\frac{|\nabla v|^2}{\alpha v} -\frac{\partial_t v}{v} 
+ \frac{\mathscr G}{v} -\frac{\beta}{\alpha} \right](x,t)
\le&~\frac{b \alpha}{t} + b \alpha \lambda_1 + 
\sqrt {b \alpha} \bigg\{\lambda_2^{4/3} + \lambda_3 + \lambda_4^2\bigg\}^{1/2}.
\end{align}
Here $M = \sup v$ over $\mathscr M \times [0, T]$ and $\lambda_j \ge 0$ $($$1 \le j \le 4$$)$ 
are as in \eqref{V-A-1.n}--\eqref{V-A-4.n}with $\mathsf K= 2c_1 \underline k_1$ in $\lambda_1$ 
and all the supremums being taken over $\mathscr M \times [0, T]$.
\end{corollary}

\begin{corollary} \label{corollary-Harnack-global.n}
Under the assumptions of Corollary $\ref{corollary-global.n}$, for every constant $\alpha >1$, 
and for all $x_1, x_2 \in \mathscr M$ and $0 <t_1<t_2<T$ we have 
\begin{align}\label{Eq5.62Har.n}
\frac{v(x_1, t_1)}{v(x_2, t_2)} \le 
{\rm exp} \left[ \frac{\alpha \mathsf L(x_1,x_2)}{4 \underline M (t_2-t_1)}
+ \frac{\mathsf H}{\alpha} (t_2 -t_1)\right] \left(\frac{t_2}{t_1}\right)^{b \alpha}.
\end{align}
Here $\underline M = \inf v$ over $\mathscr M \times [0, T]$, 
$\mathsf L=\mathsf L(x_1,x_2)$ is as in \eqref{L-Harnack-global-equation} 
and $\mathsf H$ is given by  
\begin{align}\label{V-A-H.n}
\mathsf H=&~ \lambda_0 +b \alpha^2 \lambda_1
+ \sqrt {b\alpha^3} \bigg\{\lambda_2^{4/3} + \lambda_3+ \lambda_4^2\bigg\}^{1/2},
\end{align}
with $\lambda_0 = \sup [\beta -\alpha \mathscr G/v]$ over ${\mathscr M} \times [0, T]$.
The constants in \eqref{V-A-H} are as in \eqref{estimate-2-global}.
\end{corollary}

\begin{remark} 
The following special cases and implications of the estimates in Theorem \ref{thm5.1} 
and Corollary \ref{corollary-A-Harnack} are worthy of being specifically highlighted.
\begin{itemize}
\item When $\mathscr G = \mathscr G(v)$ the terms ${\mathscr G}_x$, ${\mathscr G}_{xv}$ 
and $\Delta_\phi {\mathscr G}^x$ in $\lambda_2(R)$ and $\lambda_3(R)$ are zero. Moreover, if 
$\alpha >1$ is constant and $\beta \equiv 0$ we have $\lambda_2(R) = \sqrt{\mathsf E} \mathsf L/\alpha$, 
$\lambda_3(R) = \tilde{\mathsf M}$ and $\lambda_1(R), \lambda_4(R)$ simplify to 
\begin{align}\label{EQ-2.16}
\lambda_1(R) &= \sup_{Q_{R,T}} [\mathscr G_v+ \mathsf K]_+,\\
\lambda_4(R) &= \sqrt {\tilde{\mathsf F}}\sup_{Q_{R,T}} [(\alpha -1) (\mathscr G/v -\mathscr G_v)
- \alpha(p-1) v\mathscr G_{vv} + \tilde{\mathsf N}]_+.
\end{align}
\item In the context described above if $\mathscr G_v \le 0$ and 
$\alpha(p-1) v\mathscr G_{vv}- (\alpha -1) (\mathscr G/v -\mathscr G_v) \ge 0$ 
then $\lambda_1(R) \le \mathsf K$ and $\lambda_4(R) \le \sqrt{\tilde{\mathsf F}} \tilde{\mathsf N}$. 
In this case the local estimate \eqref{estimate-2-local} gives
\begin{align}\label{EQ-eq-2.39.n}
\frac{|\nabla v|^2}{\alpha v} -\frac{\partial_t v}{v} + \frac{\mathscr G(v)}{v}
\le&~ \frac{b \alpha}{t} + b \alpha \mathsf K 
+ \sqrt {b \alpha} \bigg\{(\sqrt{\mathsf E} \mathsf L/\alpha)^{4/3} 
+ \tilde{\mathsf M} + \tilde {\mathsf F} \tilde {\mathsf N}^2 \bigg\}^{1/2}\\
&+b(p-1) \alpha \frac{M(2R)}{R^2}[c_2 +(m-1) c_1(1+R \sqrt{k})+ 2c_1^2 ], \nonumber 
\end{align}
and the global estimate \eqref{estimate-2-global} gives 
\begin{align}\label{EQ-eq-2.40-n}
\frac{|\nabla v|^2}{\alpha v} -\frac{\partial_t v}{v} + \frac{\mathscr G(v)}{v}
\le&~ b \alpha \left\{\frac{1}{t} + 2c_1 \underline k_1\right\}\nonumber\\
&+ \sqrt {b \alpha} \bigg\{(\sqrt{\mathsf E} \mathsf L/\alpha)^{4/3} 
+ \tilde{\mathsf M} + \tilde{\mathsf F} \tilde {\mathsf N}^2 \bigg\}^{1/2}.
\end{align}
\item When $\mathscr N \equiv 0$ (equivalently $\mathscr G \equiv 0$), $\alpha >1$ 
is a constant and $\beta \equiv 0$ we have 
\begin{align}
\lambda_1(R) = \mathsf K,\quad \lambda_2(R)=\sqrt{\mathsf E}\mathsf L/\alpha, \quad
\lambda_3(R)= \tilde{\mathsf M}, \quad \lambda_4(R)= \sqrt{\tilde{\mathsf F}} \tilde{\mathsf N}. 
\end{align}
In this case the local estimate \eqref{estimate-2-local} reduces to \eqref{EQ-eq-2.39.n}
and the global estimate \eqref{estimate-2-global} reduces to \eqref{EQ-eq-2.40-n}. 
(In either case without the term $\mathscr G/v$ on the left.) 
\end{itemize}
\end{remark}

\qquad \\
{\bf The case of static triple $({\mathscr M}, g, d\mu)$.}
In the static case where the metric and potential are time independent (that is, $\partial_t g \equiv 0$ 
and $\partial_t \phi \equiv 0$) the estimates formulated above simplify further. 
Note that in this context $Q_{R,T} = {\mathscr B}_{R} \times [0, T]$ and the explicit 
bounds \eqref{vahideh-A-2}--\eqref{vahideh-A-3} are no longer needed. Specifically, here, 
we can set $\underline{k}_1=\overline{k}_1=0$, $k_2=0$ and $\ell_2=0$. In this section we discuss 
the implications of these on the estimates further.

By a close inspection of the various constants and quantities appearing on the right-hand side of 
the main estimates, it is seen that in the static case, by incorporating the above assumptions on 
$\underline{k}_1$, $\overline{k}_1$, $k_2$ and $\ell_2$ we can deduce the following. For the 
sake of convenience here we discuss the case ${\mathscr G}={\mathscr G}(v)$.

\begin{theorem}\label{thm-static-one}
Let $(\mathscr M, g, d\mu)$ be a complete smooth metric measure space with $d\mu=e^{-\phi} dv_g$. 
Assume ${\mathscr Ric}_\phi^m(g) \ge -(m-1)kg$ in ${\mathscr B}_{2R}$ for some $n \le m < \infty$ and 
$k \ge 0$. Let $u$ be a positive solution to equation \eqref{eq11} with $p>1$ and set $v =pu^{p-1}/(p-1)$. 
Then for every $(x,t)$ in ${\mathscr B}_R \times (0, T]$ we have 
\begin{align} \label{estimate-5.4-static-one}
&\left[\frac{|\nabla v|^2}{\alpha v} -\frac{\partial_t v}{v} + \frac{\mathscr G (v)}{v} \right](x,t)
\le \frac{b \alpha}{t} + b \alpha \sup_{{\mathscr B}_{2R} \times [0, T]} \left[\mathscr G_v + \frac{\alpha'}{\alpha} 
+ \frac{b \alpha^2 p^2 M(2R)}{2(\alpha -1)} \frac{c_1^2}{R^2} \right]_+ \nonumber \\
&+b(p-1) \alpha \frac{M(2R)}{R^2} [c_2 +(m-1) c_1(1+R \sqrt{k})+ 2c_1^2]\\
&+ b  \sup_{{\mathscr B}_{2R} \times [0, T]} \left[ \frac{\alpha}{2} \left(\frac{\mathscr G}{v} -\mathscr G_v\right)
- \frac{\alpha^2(p-1)}{2(\alpha -1)} v \mathscr G_{vv}
+ \frac{2\alpha(p-1)M(2R)(m-1) k - \alpha'}{2(\alpha -1)} \right] _+. \nonumber
\end{align}
Here $M = \sup v$ over ${\mathscr B}_{2R} \times [0, T]$, $c_1$, $c_2 >0$ are the constants in \eqref{9.30}.  
\end{theorem}

Note that here by referring to \eqref{EQ-eq-K}--\eqref{EQ-eq-E} we have 
$\mathsf K = c_1^2 b \alpha^2 p^2 M /[2(\alpha -1)R^2]$, 
$\mathsf L = 0$, $\mathsf M = 0$, $\mathsf N = 2(p-1)M(m-1) k$ 
and $\mathsf F = b \alpha^2 /[4 (\alpha -1)^2-2 \varepsilon b \alpha^2]$. 
Additionally, as $\beta=0$ and ${\mathscr G}={\mathscr G}(v)$, by incorporating 
these into \eqref {V-A-1}--\eqref {V-A-4} and noting $\mu_2(R)=0$, $\mu_3(R)=0$, 
the conclusion \eqref{estimate-5.4-static-one} follows at once by letting $\varepsilon \searrow 0$ in \eqref {estimate-5.4}. 
The passage to the limit $\varepsilon \searrow 0$ is allowed since here 
$\mu_2(R)=0$ (in turn as a consequence of ${\mathscr G}_x \equiv 0$ and ${\mathsf L}=0$).

Letting $R \nearrow \infty$ or alternatively using directly the global estimate 
\eqref {estimate-5.4-global} with these values leads to

\begin{corollary}
If ${\mathscr Ric}_\phi^m(g) \ge -(m-1)kg$ holds everywhere on $\mathscr M$ and $v$ is bounded 
with $M = \sup v$ over $\mathscr M \times [0, T]$ then for every $x \in \mathscr M$ and 
$0< t \le T$ we have
\begin{align}
&\left[\frac{|\nabla v|^2}{\alpha v} -\frac{\partial_t v}{v} + \frac{\mathscr G (v)}{v} \right](x,t)
\le \frac{b \alpha}{t} + b \alpha \sup_{{\mathscr M} \times [0, T]} \left[\mathscr G_v 
+ \frac{\alpha'}{\alpha} \right]_+ \nonumber \\
&+ b  \sup_{{\mathscr M} \times [0,T]} \left[ \frac{\alpha}{2} \left(\frac{\mathscr G}{v} -\mathscr G_v\right)
- \frac{\alpha^2(p-1)}{2(\alpha -1)} v \mathscr G_{vv}
+ \frac{2\alpha(p-1)M(m-1) k - \alpha'}{2(\alpha -1)} \right] _+. 
\end{align} 
\end{corollary}

Similar to what has been described above we have the following counterparts of the second main 
estimate, i.e., Theorem \ref{Thrm-2.5} in the static case.

\begin{theorem}\label{thm-static-two}
Let $(\mathscr M, g, d\mu)$ be a complete smooth metric measure space with $d\mu=e^{-\phi} dv_g$. 
Assume ${\mathscr Ric}_\phi^m(g) \ge -(m-1)kg$ in ${\mathscr B}_{2R}$ for some $n \le m < \infty$ and 
$k \ge 0$. Let $u$ be a positive solution to equation \eqref{eq11} with $p>1$ and set $v =pu^{p-1}/(p-1)$. 
Then for every $(x,t)$ in ${\mathscr B}_R \times (0, T]$ we have 
\begin{align} \label{estimate-5.4-static-two}
&\left[\frac{|\nabla v|^2}{\alpha v} -\frac{\partial_t v}{v} + \frac{\mathscr G(v)}{v} \right](x,t)
\le \frac{b \alpha}{t} + b \alpha \sup_{{\mathscr B}_{2R} \times [0, T]} 
\left[\mathscr G_v + \frac{b \alpha^2 p^2 M(2R)}{2(\alpha -1)} \frac{c_1^2}{R^2} \right]_+ \nonumber \\ 
&+ b (p-1)\alpha \frac{M(2R)}{R^2} [c_2 +(m-1) c_1(1+R \sqrt{k})+ 2c_1^2] + b \sqrt{\alpha} \times \\
& \sup_{{\mathscr B}_{2R} \times [0, T]} \left[ \frac{\sqrt \alpha}{2} \left(\frac{\mathscr G}{v} -\mathscr G_v\right)
- \frac{\alpha^\frac{3}{2}(p-1)}{2(\alpha -1)} v \mathscr G_{vv}
+ \frac{2\alpha(p-1)M(2R)(m-1) k - \alpha'}{2\sqrt{\alpha} (\alpha -1)} \right] _+. \nonumber
\end{align}
Here $M = \sup v$ over ${\mathscr B}_{2R} \times [0, T]$, $c_1$, $c_2 >0$ are the constants in \eqref{9.30}.  
\end{theorem}

Note that similar to Theorem \ref{thm-static-one}, here, by referring to \eqref{EQ-eq-K}--\eqref{EQ-eq-L} 
and \eqref{EQ-eq-M.n}--\eqref{EQ-eq-E.n} we have $\mathsf K = c_1^2 b \alpha^2 p^2 M /[2(\alpha -1)R^2]$, 
$\mathsf L =0$, $\tilde {\mathsf M} = 0$, $\tilde{\mathsf N} = 2(p-1)M(m-1) k/\alpha$, 
and $\tilde{\mathsf F} = b \alpha^3 /[4 (\alpha -1)^2-2 \varepsilon b \alpha^3]$. Since 
$\beta=0$ and ${\mathscr G}={\mathscr G}(v)$, by incorporating these into 
\eqref {V-A-1.n}--\eqref {V-A-4.n} and noting $\lambda_2(R)=0$, $\lambda_3(R)=0$, 
the conclusion \eqref{estimate-5.4-static-two} follows by letting $\varepsilon \searrow 0$ 
in \eqref{estimate-2-local}. 
The passage to the limit $\varepsilon \searrow 0$ is allowed as 
$\lambda_2(R)=0$ (in turn as a consequence of ${\mathscr G}_x \equiv 0$ and ${\mathsf L}=0$).

Letting $R \nearrow \infty$ or alternatively using directly the global estimate 
\eqref{estimate-2-global} with these values leads to

\begin{corollary}
If ${\mathscr Ric}_\phi^m(g) \ge -(m-1)kg$ holds everywhere on $\mathscr M$ and $v$ is bounded 
with $M = \sup v$ over $\mathscr M \times [0, T]$ then for every $x \in \mathscr M$ and 
$0< t \le T$ we have
\begin{align} \label{global-static-2.n}
&\left[\frac{|\nabla v|^2}{\alpha v} -\frac{\partial_t v}{v} + \frac{\mathscr G (v)}{v} \right](x,t)
\le \frac{b \alpha}{t} + b \alpha \sup_{{\mathscr M} \times [0, T]} 
[\mathscr G_v]_+ + b \sqrt{\alpha} \times \nonumber \\ 
&\times \sup_{{\mathscr M} \times [0, T]} \left[ \frac{\sqrt \alpha}{2} \left(\frac{\mathscr G}{v} -\mathscr G_v\right)
- \frac{\alpha^\frac{3}{2}(p-1)}{2(\alpha -1)} v \mathscr G_{vv}
+ \frac{2\alpha(p-1)M(m-1) k - \alpha'}{2\sqrt{\alpha} (\alpha -1)} \right] _+. 
\end{align}
\end{corollary}

\qquad \\
{\bf Some further remarks and comments.} Having formulated the two set of gradient estimates and their 
consequences we now discuss some special choices of the functions $\alpha=\alpha(t)$, $\beta=\beta(t)$ 
that lead to certain simplifications in the estimates. Note that here we write $\gamma=(m-1)k(p-1)M$ for 
brevity. It is also interesting to compare these with the special cases 
in \cite{Ha93, Huang-Huang-Li, Hunag-Li, Huang}.

\begin{remark} \label{interesting-alpha-beta-remark} Some examples of time dependent Harnack coefficient functions 
$\alpha=\alpha(t)$ and $\beta=\beta(t)$ leading to certain simplifications in the estimates are given below: 
\begin{itemize}
\item $\alpha(t)= e^{2\gamma t}$ [and any $\beta(t)$] gives $2\gamma \alpha-\alpha'=0$.
\item $\alpha(t) = 1 + [\cosh (\gamma t) \sinh (\gamma t) - \gamma t]/\sinh^2 (\gamma t)$ and 
$\beta(t) = b \gamma [\coth (\gamma t) + 1]$. Then $\alpha$ and $\beta$ satisfy the following equations
\begin{align}
2 \beta/b-\alpha' = 2\alpha (\beta/b -\gamma), \qquad \beta' + 2(\beta/b-\gamma) \beta - \beta^2/b=0. 
\end{align} 
 
\item $\alpha(t) = 1+ 2 \gamma t/3$ and $\beta(t) = b(1/t + \gamma + \gamma^2 t/3)$. 
Then $\alpha$ and $\beta$ satisfy the following equations
\begin{align}
2(1/t + \gamma)-\alpha' = 2\alpha/t, \qquad \beta' + 2 \beta/t - b (1/t + \gamma)^2 =0. 
\end{align} 
\end{itemize}
\end{remark}
The first choice of $\alpha$ removes the $k$ dependent terms in $\lambda_4(R)$ [{\it see} \eqref{V-A-4.n}] and hence simplifies 
\eqref{estimate-2-local}, \eqref {estimate-2-global} and \eqref{estimate-5.4-static-two}. It also eliminates the 
curvature term in \eqref{global-static-2.n} and so subject to  
$\mathscr G_v \le 0$ and $\alpha(p-1) v\mathscr G_{vv}- (\alpha -1) (\mathscr G/v -\mathscr G_v) \ge 0$ 
reduces the right-hand side of the global estimate \eqref{global-static-2.n} to $b \alpha(t)/t$. This is a far 
reaching generalisation of the results in \cite{Huang-Huang-Li, Hunag-Li} on the equation \eqref {EQ-1.6} 
to the nonlinear porous medium equation \eqref{eq11}. The other cases are similar and will be discussed in due course.

\qquad \\
{\bf A special nonlinearity.} As the case $\mathscr G=\mathscr G(v)$ along with the conditions 
$\mathscr G_v \le 0$ and $\alpha(p-1) v\mathscr G_{vv}- (\alpha -1) (\mathscr G/v -\mathscr G_v) \ge 0$ have 
proven to be important (see the remarks following the estimates) let us now give some useful examples where these 
conditions can be verified. Consider the nonlinearity 
\begin{align}
\mathscr G(v) = \sum_{j=1}^N \mathsf{A}_j v^{a_j} + \sum_{j=1}^N \mathsf{B}_j v^{b_j}, \qquad v>0,
\end{align}
where $\mathsf A_j \ge0$ , $\mathsf B_j \le 0$ and $a_j, b_j$ are real exponents. 
Then a straightforward calculation gives 
\begin{align}
\mathscr G_v (v) =\sum_{j=1}^N \mathsf{A}_j a_j v^{a_j-1} + \sum_{j=1}^N \mathsf{B}_j b_j v^{b_j-1},
\end{align}
and so $\mathscr G_v \le 0$ when $a_j \le0$ and $b_j \ge 0$. Moreover a further differentiation gives
\begin{align}
\alpha(p-1) v^2\mathscr G_{vv}(v)&- (\alpha -1) [\mathscr G(v) -v\mathscr G_v(v)]\nonumber\\
= &~\alpha(p-1) \left(\sum_{j=1}^N \mathsf{A}_j a_j (a_j-1)v^{a_j} + \sum_{j=1}^N \mathsf{B}_j b_j (b_j-1)v^{b_j}\right)\nonumber\\
& + (\alpha -1) \left(\sum_{j=1}^N \mathsf{A}_j (a_j-1)v^{a_j} + \sum_{j=1}^N \mathsf{B}_j (b_j-1)v^{b_j}\right)\nonumber\\
=&~ \sum_{j=1}^N  \mathsf{A}_j (a_j-1) [\alpha(p-1) a_j+ (\alpha-1)] v^{a_j}\nonumber\\
& +\sum_{j=1}^N  \mathsf{B}_j (b_j-1) [\alpha(p-1) b_j+ (\alpha-1)] v^{b_j}.
\end{align}
Thus based on the above ranges of $a_j$, $b_j$ we have 
$\alpha(p-1) v\mathscr G_{vv}- (\alpha -1) (\mathscr G/v -\mathscr G_v) \ge 0$ 
when $a_j \le (1-\alpha)/[\alpha(p-1)]$ and $0\le b_j \le 1$. Recall that $\alpha>1$, $p>1$.

\section{Evolution lemmas and the Harnack quantity $F$} 
\label{sec3}

In this section we introduce a Harnack quantity that is built out of the positive solution $u= u(x,t)$ and establish a parabolic inequality 
by considering its evolution under the operator $\mathscr L^p_v$ defined by 
\begin{equation}\label{eq5.6}
\mathscr L^p_v = \partial _t -(p-1) v \Delta_ \phi.
\end{equation}
Here as before $p>1$ and $v(x,t) =pu^{p-1}(x,t)/(p-1)$.

\begin{lemma}\label{Lemma-5.1}
Let $u$ be a positive solution to the equation \eqref{eq11}.
Then the function $v$ satisfies the equation 
\begin{align}\label{eq-8.1}
\mathscr L^p_v[v] =[\partial_t - (p-1) v \Delta_\phi] v = |\nabla v|^2+\mathscr G(t,x,v),
\end{align}
where $\mathscr G=\mathscr G(t,x,v)$ is as in \eqref{eq5.2}. 
\end{lemma}

\begin {proof}
As here $u = [(p-1)v/p]^{1/(p-1)}$ a  straightforward calculation gives
\begin{align}
\partial_t u = (1/p) [(p-1)v/p]^{\frac{2-p}{p-1}} \partial_t v, \qquad 
\nabla u^p = [(p-1)v/p]^{\frac{1}{p-1}} \nabla v, 
\end{align}  
\begin{align}
\Delta u^p = (1/p) [(p-1)v/p]^{\frac{2-p}{p-1}}[(p-1) v \Delta v + |\nabla v|^2].
\end{align}
Therefore by recalling the relation 
$\Delta_\phi w = \Delta w - \langle \nabla \phi, \nabla w \rangle$ from \eqref{f-Lap-definition} it is seen that 
\begin{align}
\Delta_\phi u^p 
&= (1/p) [(p-1)v/p]^{\frac{2-p}{p-1}}[(p-1) v \Delta v + |\nabla v|^2]
- [(p-1)v/p]^{\frac{1}{p-1}}\langle \nabla \phi , \nabla v \rangle\nonumber\\
&= (1/p) [(p-1)v/p]^{\frac{2-p}{p-1}}[(p-1) v \Delta_\phi v + |\nabla v|^2].
\end{align}
Substituting the above in \eqref{eq11} and rearranging terms gives the desired conclusion.  
\end{proof}

\begin{lemma} \label{Lpv-quotient-lemma}
For a pair of functions $f$, $g$ of class $\mathscr C^2(M)$ with $g$ nowhere zero we have the identity 
\begin{align}\label{eq5.19}
\mathscr L^p_v [f/g] = 
(1/g) \mathscr L^p_v [f] + 2 (p-1) v \langle \nabla (f/g), \nabla g/g \rangle - (f/g^2) \mathscr L^p_v [g].
\end{align}
Note that when $g$ is positive we can write $\langle \nabla (f/g), \nabla g/g \rangle= \langle \nabla (f/g), \nabla \log g \rangle$.
\end{lemma}

\begin{proof}
Referring to \eqref{eq5.6} we can write $\mathscr L^p_v [f/g] = [\partial _t - (p-1) v \Delta_ \phi ] (f/g)$ and so making use 
of $\Delta_\phi (f/g) =(1/g) \Delta_\phi f - 2 (1/g^2) [\langle \nabla f, \nabla g \rangle - (f/g)|\nabla g|^2] - (f/g^2) \Delta_\phi g$ we have  
\begin{align} 
\mathscr L^p_v [f/g] 
=&~ (1/g) \partial_t f - (f/g^2) \partial_t g - (p-1)v (1/g) \Delta_ \phi f \nonumber\\
& +2(p-1) v(1/g^2) [\langle \nabla f, \nabla g \rangle - (f/g) |\nabla g|^2]+ (p-1)v (f/g^2) \Delta_\phi g\nonumber\\
=&~(1/g) \mathscr L^p_v [f] +2(p-1) v(1/g^2) [\langle \nabla f, \nabla g \rangle - (f/g) |\nabla g|^2]
- (f/g^2) \mathscr L^p_v [g]\nonumber\\
=&~ (1/g) \mathscr L^p_v [f] +2 (p-1) v \langle \nabla (f/g), \nabla g/g \rangle
- (f/g^2) \mathscr L^p_v [g]. 
\end{align} 
The last identity makes use of $(1/g^2)[\langle \nabla f, \nabla g \rangle -(f/g)|\nabla g|^2] 
= \langle \nabla (f/g), \nabla g/g \rangle$.
\end{proof}

\begin{lemma} \label{geometric-evolution-lemma-two}
For every pair of smooth functions $\phi=\phi(x,t)$ and $v=v(x,t)$ we have 
\begin {align} \label{Lap-f-evolve-equation}
\partial_t \Delta_\phi v - \Delta_\phi \partial_t v 
=&~2 \langle h , \nabla \nabla v \rangle - \langle 2 {\rm div} h - \nabla ({\rm Tr}_g h),\nabla v \rangle \nonumber \\
&~+ 2h (\nabla \phi , \nabla v) - \langle \nabla \partial_t \phi , \nabla v \rangle. 
\end{align}
\end{lemma}
\begin{proof}
See the proofs of Lemma~5.1 and Lemma~5.2 (on pp.~21-22) in \cite{TVahNA}. 
\end{proof}

\begin{lemma} \label{lemma-5.6-Equality}
Let $u=u(x,t)$ be a positive solution to equation \eqref{eq11} and let $F=F(x,t)$ be defined by 
\begin{align}\label{eq5.21}
 F = \frac{|\nabla v|^2}{v} -\alpha(t) \frac{\partial_t v}{v} + \alpha(t) \frac{\mathscr G(t,x,v)}{v} -\beta(t),
 \end{align}
where $v =p u^{p-1}/(p-1)$, $\mathscr G(t,x,v)$ is as in \eqref{eq5.2} and $\alpha=\alpha(t)$, 
$\beta=\beta(t)$ are arbitrary functions of class $\mathscr{C}^1$. Suppose that the metric and 
potential are time dependent, of class $\mathscr{C}^2$ and let $h$ be the tensor field in $\eqref{2.4}$. Then
\begin{align}\label{eq-5.27}
\mathscr L^p_v [F] =& - [(p-1) \Delta_\phi v]^2 
+ 2p \langle \nabla F, \nabla v \rangle
- 2(p-1) |\nabla \nabla v|^2
-2(p-1) {\mathscr Ric}_\phi^m (\nabla v, \nabla v) \nonumber\\
&-2(p-1)\frac{\langle \nabla \phi , \nabla v \rangle^2}{m-n}
+\frac{2}{v}[\alpha(t) -1] h ( \nabla v, \nabla v)
+2\alpha(t) (p-1) \langle h , \nabla \nabla v \rangle \nonumber\\
&+ \alpha(t) (p-1) \langle 2 {\rm div} h - \nabla ({\rm Tr}_g h),\nabla v \rangle
+ \alpha(t) (p-1) \langle \nabla \partial_t \phi , \nabla v \rangle-\beta' (t)\nonumber\\
&- 2 \alpha(t) (p-1)h (\nabla \phi , \nabla v)
-\alpha'(t) \frac{\partial_t v}{v}
 + [1-\alpha (t)] \left[ \frac{\partial_t v}{v} - \frac{\mathscr G(t,x,v)}{v} \right]^2 \nonumber\\
&+\frac{2}{v}[1-\alpha(t)] \langle \nabla \mathscr G(t,x,v), \nabla v \rangle
- \alpha (t)(p-1) \Delta_\phi \mathscr G(t,x,v) \nonumber\\
&+[\alpha(t)-1] \frac{|\nabla v|^2}{v} \frac{\mathscr G(t,x,v)}{v}+ \alpha'(t) \frac{\mathscr G(t,x,v)}{v}.
\end{align}
\end{lemma}

\begin{proof}
Using linearity and the explicit formulation of $F$ in \eqref{eq5.21} we can write
\begin{align}\label{eq-2.7}
\mathscr L ^p_v [F] =&~ \mathscr L ^p_v [|\nabla v|^2/v -\alpha(t) \partial_t v/v + \alpha(t) \mathscr G(t,x,v)/v -\beta(t)]\nonumber\\
=&~ \mathscr L^p_v [|\nabla v|^2/v ] - \alpha (t) \mathscr L^p_v [\partial_t v/v]
+ \alpha(t) \mathscr L^p_v [\mathscr G(t,x,v)/v]\nonumber\\
&-\alpha'(t) \partial_t v/v
+ \alpha'(t) \mathscr G(t,x,v)/v -\beta' (t).
\end{align}

Now we expand the expression on the right-hand side by working out separately the first three terms involving 
$\mathscr L ^p_v$ by invoking Lemma \ref{Lpv-quotient-lemma}. Using the weighted Bochner-Weitzenb\"ock formula \eqref{Bochner-1}
we can write 
\begin{align}\label{eq2.13}
\mathscr L^p_v [|\nabla v|^2] =&~ [\partial_t -(p-1) v \Delta_\phi] |\nabla v|^2 \nonumber\\
= & -[\partial_t g] ( \nabla v, \nabla v) +2 \langle \nabla v ,\nabla \partial_t v \rangle \nonumber\\
&-2 (p-1) v [ |\nabla \nabla v|^2 + \langle \nabla v ,\nabla \Delta_\phi v \rangle 
+ {\mathscr Ric}_\phi^m(\nabla v, \nabla v) +\langle \nabla \phi , \nabla v \rangle^2/(m-n)] \nonumber\\
= & -2h ( \nabla v, \nabla v) +2 \left \langle \nabla v, \nabla [(p-1) v \Delta_\phi v 
+ |\nabla v|^2 +\mathscr G(t,x,v) ] \right\rangle \nonumber\\
& -2 (p-1) v [ |\nabla \nabla v|^2 + \langle \nabla v, \nabla \Delta_\phi v \rangle
+ {\mathscr Ric}_\phi^m (\nabla v, \nabla v) + \langle \nabla \phi , \nabla v \rangle^2/(m-n)] \nonumber\\
= & -2h ( \nabla v, \nabla v) +2 (p-1) |\nabla v|^2 \Delta_\phi v
 + 2 (p-1) v \langle \nabla v ,\nabla \Delta_\phi v \rangle\nonumber\\
&+ 2\langle \nabla v, \nabla |\nabla v|^2 \rangle +2 \langle \nabla v ,\nabla \mathscr G(t,x,v) \rangle
- 2(p-1) v |\nabla \nabla v|^2 \nonumber\\
&- 2 (p-1) v \langle \nabla v, \nabla \Delta_\phi v \rangle
-2(p-1)v {\mathscr Ric}_\phi^m (\nabla v, \nabla v) \nonumber\\
&- 2 (p-1) v \langle \nabla \phi , \nabla v \rangle^2/(m-n)\nonumber\\
= & -2h ( \nabla v, \nabla v) +2 (p-1) |\nabla v|^2 \Delta_\phi v + 2\langle \nabla v, \nabla |\nabla v|^2 \rangle\nonumber\\
&+2 \langle \nabla v, \nabla \mathscr G(t,x,v) \rangle-2(p-1) v |\nabla \nabla v|^2\nonumber\\
& -2(p-1)v {\mathscr Ric}_\phi^m (\nabla v, \nabla v) - 2 (p-1) v \langle \nabla \phi , \nabla v \rangle^2/(m-n).
\end{align}

Hence, for the first term on the right-hand side of \eqref{eq-2.7}, upon making note of Lemma \ref{Lpv-quotient-lemma} 
and substituting from \eqref{eq-8.1} and \eqref{eq2.13}, we have
\begin{align}\label{eq-2.10}
\mathscr L^p_v [|\nabla v|^2/v] =&~ (1/v) \mathscr L^p_v [|\nabla v|^2] 
+2(p-1)v \langle \nabla (|\nabla v|^2/v ), \nabla v/v \rangle 
- (|\nabla v|^2/v^2) \mathscr L^p_v [v] \nonumber\\
=& - 2 h ( \nabla v, \nabla v)/v+2 (p-1) |\nabla v|^2 \Delta_\phi v/v
+2 \langle \nabla |\nabla v|^2, \nabla v \rangle/v\nonumber\\
&+2 \langle \nabla v ,\nabla \mathscr G(t,x,v) \rangle/v
-2(p-1) |\nabla \nabla v|^2 -2(p-1) {\mathscr Ric}_\phi^m (\nabla v, \nabla v) \nonumber\\
&- 2 (p-1) \langle \nabla \phi , \nabla v \rangle^2/(m-n)
+2(p-1)v \langle \nabla (|\nabla v|^2/v), \nabla v/v \rangle \nonumber\\
&- |\nabla v|^4/v^2 - (|\nabla v|^2/v^2) \mathscr G(t,x,v).
\end{align}
Similarly, to deal with the  second term on the right-hand side of \eqref{eq-2.7} 
by tacking advantage of \eqref{Lap-f-evolve-equation} 
in Lemma \ref{geometric-evolution-lemma-two}, we start with
\begin{align}\label{eq2.16}
\mathscr L^p_v [\partial_t v] =&~ [\partial_t -(p-1) v \Delta_\phi] \partial_t v \nonumber\\
=&~ \partial_t [(p-1) v \Delta_\phi v + |\nabla v|^2 + \mathscr G (t,x,v) ]-(p-1) v \Delta_\phi \partial_t v \nonumber\\
=&~ (p-1) \partial_t v \Delta_\phi v +(p-1)v \partial_t \Delta_\phi v -2h ( \nabla v, \nabla v)\nonumber\\
&+ 2 \langle \nabla v  ,\nabla \partial_t v \rangle 
+\partial_t \mathscr G (t,x,v) -(p-1) v \Delta_\phi \partial_t v\nonumber\\
=& ~(p-1) \partial_t v \Delta_\phi v -2(p-1)v \langle h , \nabla \nabla v \rangle 
- (p-1)v \langle 2 {\rm div} h - \nabla ({\rm Tr}_g h),\nabla v \rangle \nonumber\\
&- (p-1)v \langle \nabla v, \nabla \partial_t \phi \rangle + 2(p-1)v h(\nabla \phi , \nabla v)
-2h ( \nabla v, \nabla v)\nonumber\\
& + 2 \langle \nabla v, \nabla \partial_t v \rangle 
+\partial_t \mathscr G (t,x,v).
\end{align}
Next, by using Lemma \ref{Lpv-quotient-lemma} and substituting from \eqref {eq-8.1} and \eqref{eq2.16} we have 
\begin{align}\label{eq-2.8}
\mathscr L^p_v [\partial_t v/v] =&~ (1/v) \mathscr L^p_v [\partial_t v] 
+2(p-1)v \langle \nabla (\partial_t v/v), \nabla v/v \rangle
- (\partial_t v/v^2) \mathscr L^p_v [v] \nonumber\\
=&~ (p-1) \partial_t v \Delta_\phi v/v  -2(p-1) \langle h , \nabla \nabla v \rangle 
- (p-1) \langle 2 {\rm div} h - \nabla ({\rm Tr}_g h),\nabla v \rangle \nonumber\\
& - (p-1) \langle \nabla v, \nabla \partial_t \phi \rangle + 2(p-1)h (\nabla \phi , \nabla v)
- 2 h ( \nabla v, \nabla v)/v\nonumber\\
&+ 2 \langle \nabla v , \nabla \partial_t v\rangle /v
+ \partial_t \mathscr G (t,x,v) /v +2(p-1)v \langle \nabla (\partial_t v/v), \nabla v/v \rangle\nonumber\\
& - (\partial_t v/v^2) |\nabla v|^2 - (\partial_t v/v^2) \mathscr G(t,x,v) .
\end{align}
Finally for the third term on the right-hand side of \eqref{eq-2.7}, again by using 
Lemma \ref{Lpv-quotient-lemma} and substituting from \eqref {eq-8.1}, we can write 
\begin{align}\label{eq-2.13}
\mathscr L^p_v [\mathscr G(t,x,v)/v]=&~ (1/v) \mathscr L^p_v [\mathscr G(t,x,v)] 
-(\mathscr G(t,x,v)/v^2) \mathscr L^p_v [v] \nonumber\\
&+2(p-1)v \langle \nabla (\mathscr G(t,x,v)/v), \nabla v/v \rangle\nonumber\\
=&~ \partial_t \mathscr G(t,x,v)/v -(p-1) \Delta_\phi \mathscr G(t,x,v)\nonumber\\
&+2(p-1) v \langle \nabla (\mathscr G(t,x,v)/v), \nabla v/v \rangle \nonumber\\
&- (\mathscr G(t,x,v)/v^2) |\nabla v|^2 - [\mathscr G(t,x,v)/v]^2.
\end{align}
Substituting \eqref{eq-2.10}, \eqref{eq-2.8} and \eqref{eq-2.13} back into \eqref{eq-2.7} and rearranging terms gives 
\begin{align}\label{eq-2.12}
\mathscr L^p_v [F] =&~ \mathscr L ^p_v [|\nabla v|^2/v] 
-\alpha(t) \mathscr L^p_v [\partial_t v/v]
+ \alpha(t) \mathscr L^p_v [\mathscr G(t,x,v)/v]\\
&-\alpha'(t) \partial_t v/v
+ \alpha'(t) \mathscr G(t,x,v)/v -\beta' (t) \nonumber\\
=& - 2 h ( \nabla v, \nabla v)/v+2 (p-1) |\nabla v|^2 \Delta_\phi v /v
+2 \langle \nabla |\nabla v|^2, \nabla v \rangle/v\nonumber\\
&+ 2 \langle \nabla v ,\nabla \mathscr G(t,x,v) \rangle/v - 2(p-1) |\nabla \nabla v|^2 
-2(p-1) {\mathscr Ric}_\phi^m (\nabla v, \nabla v)\nonumber\\
&- 2 (p-1) \langle \nabla \phi , \nabla v \rangle^2/(m-n)
+2(p-1) v \langle \nabla (|\nabla v|^2/v), \nabla v/v  \rangle \nonumber\\
& - |\nabla v|^4/v^2 - (|\nabla v|^2/v^2) \mathscr G(t,x,v)
- \alpha(t) (p-1) \partial_t v \Delta_\phi v/v\nonumber\\
&+2\alpha(t) (p-1) \langle h , \nabla \nabla v \rangle 
+ \alpha(t) (p-1) \langle 2 {\rm div} h - \nabla ({\rm Tr}_g h),\nabla v \rangle\nonumber\\
&+ \alpha(t) (p-1) \langle \nabla v, \nabla \partial_t \phi \rangle - 2 \alpha(t) (p-1)h (\nabla \phi , \nabla v)
+ 2\alpha(t) h ( \nabla v, \nabla v)/v\nonumber\\
&- 2\alpha(t) \langle \nabla v, \nabla \partial_t v \rangle /v
-\alpha (t) \partial_t \mathscr G (t,x,v)/v
-2 \alpha(t)(p-1) v \langle \nabla ( \partial_t v/v), \nabla v/v \rangle\nonumber\\
&+ \alpha(t) (\partial_t v/v^2) |\nabla v|^2 + \alpha(t) (\partial_t v/v^2) \mathscr G(t,x,v)  
+ \alpha(t) \partial_t \mathscr G(t,x,v)/v \nonumber\\
&- \alpha(t)(p-1) \Delta_\phi \mathscr G(t,x,v) 
+2 \alpha(t)(p-1) v \langle \nabla (\mathscr G(t,x,v)/v), \nabla v/v \rangle\nonumber\\
&-\alpha(t) (|\nabla v|^2/v^2) \mathscr G(t,x,v)- \alpha(t) [\mathscr G(t,x,v)/v]^2\nonumber\\
&-\alpha'(t) \partial_t v/v+ \alpha'(t) \mathscr G(t,x,v)/v -\beta' (t).\nonumber
\end{align}

Taking into account the cancellation of terms and a further rearrangement of the identity then results in 
\begin{align}\label{eq-2.20}
\mathscr L ^p_v [F] =&~ [2 |\nabla v|^2/v 
- \alpha(t) \partial_t v/v] (p-1) \Delta_\phi v 
- |\nabla v|^4/v^2 + \alpha(t) \partial_t v |\nabla v|^2/v^2\nonumber\\
& -2(p-1) |\nabla \nabla v|^2 + 2 [\alpha(t) -1] h ( \nabla v, \nabla v)/v
+ 2 \langle \nabla v ,\nabla \mathscr G(t,x,v) \rangle /v\nonumber\\
&-2(p-1) {\mathscr Ric}_\phi^m (\nabla v, \nabla v) -2 (p-1) \langle \nabla \phi , \nabla v \rangle^2/(m-n) \nonumber\\
&+2\alpha(t) (p-1) \langle h , \nabla \nabla v \rangle-[1+\alpha(t)] |\nabla v|^2 \mathscr G(t,x,v)/v^2\nonumber\\
&+2(p-1) v\langle \nabla (|\nabla v|^2/v), \nabla v/v \rangle + \alpha(t) (p-1) \langle \nabla \partial_t \phi , \nabla v \rangle\nonumber\\
&+ 2 \langle \nabla v, \nabla |\nabla v|^2 \rangle/v- 2 \alpha(t) (p-1)h (\nabla \phi , \nabla v)
- 2\alpha(t) \langle \nabla \partial_t v, \nabla v \rangle/v\nonumber\\
&+ \alpha(t) (p-1) \langle 2 {\rm div} h - \nabla ({\rm Tr}_g h),\nabla v \rangle + \alpha(t) \partial_t v \mathscr G(t,x,v)/v^2\nonumber\\
&-2 \alpha(t)(p-1) v \langle \nabla (\partial_t v/v), \nabla v/v \rangle- \alpha (t)(p-1) \Delta_\phi \mathscr G(t,x,v)\nonumber\\
&+2\alpha(t)(p-1)v \langle \nabla (\mathscr G(t,x,v)/v), \nabla v/v \rangle-\alpha (t) [\mathscr G(t,x,v)/v]^2\nonumber\\
&-\alpha'(t) \partial_t v/v+ \alpha'(t) \mathscr G(t,x,v)/v -\beta' (t).
\end{align}

Recalling that $F = |\nabla v|^2/v -\alpha(t) \partial_t v/v + \alpha(t) \mathscr G(t,x,v)/v -\beta(t)$, 
a close inspection of the terms on the right-hand side, enable us to make use of  
\begin{align}\label{eq-2.15}
2(p-1)v \{&\langle \nabla(|\nabla v|^2/v), \nabla v/v \rangle
-\alpha(t) \langle\nabla ( \partial_t v/v), \nabla v/v \rangle \nonumber\\
&+\alpha(t)  \langle \nabla (\mathscr G(t,x,v)/v), \nabla v/v \rangle \}
= 2(p-1) \langle \nabla F , \nabla v \rangle,
\end{align}
and 
\begin{align}\label{eq-2.22}
2 \langle \nabla |\nabla v|^2, \nabla v \rangle/v
-2 \alpha(t) \langle \nabla \partial_t v, \nabla v \rangle/v
=&~ 2 \langle \nabla [(F+ \beta(t))v], \nabla v \rangle/v\nonumber\\
&- 2\alpha (t) \langle \nabla \mathscr G(t,x,v), \nabla v \rangle/v\nonumber\\
=&~ 2 [F+ \beta(t)] |\nabla v|^2/v + 2 \langle \nabla F, \nabla v \rangle\nonumber\\
&-2 \alpha (t) \langle \nabla \mathscr G(t,x,v), \nabla v\rangle/v.
\end{align}
Putting \eqref{eq-2.15}--\eqref{eq-2.22} together gives 
\begin{align}\label{eq-2.17}
&2(p-1)v \{ \langle \nabla (|\nabla v|^2/v), \nabla v/v \rangle
- \alpha(t) \langle\nabla ( \partial_t v/v), \nabla v/v \rangle \nonumber\\
&+\alpha(t) \langle \nabla (\mathscr G(t,x,v)/v), \nabla v/v \rangle \}
+ 2\langle \nabla v, \nabla |\nabla v|^2 \rangle/v
- 2\alpha(t) \langle \nabla v , \nabla \partial_t v \rangle/v\nonumber\\
&= 2p \langle \nabla F, \nabla v \rangle + 2 [F+ \beta(t)] |\nabla v|^2/v
- 2\alpha (t) \langle \nabla v ,\nabla \mathscr G(t,x,v) \rangle/v\nonumber\\
&= 2p \langle \nabla F, \nabla v \rangle 
+2 \left[|\nabla v|^2/v -\alpha(t) \partial_t v/v
+ \alpha(t) \mathscr G(t,x,v)/v \right] |\nabla v|^2/v\nonumber\\
&-2\alpha (t) \langle \nabla v ,\nabla \mathscr G(t,x,v) \rangle/v.
\end{align}
Using \eqref{eq-8.1} for the first three terms on the right-hand side of \eqref{eq-2.20} we can write 
\begin{align}\label{eq-2.18}
[ 2 |\nabla v|^2/v - \alpha(t) \partial_t v/v]& (p-1) \Delta_\phi v -
|\nabla v|^4/v^2 + \alpha(t)\partial_t v |\nabla v|^2/v^2 \nonumber\\
=&~ [2 |\nabla v|^2/v - \alpha(t) \partial_t v/v] 
[\partial_t v/v -|\nabla v|^2/v - \mathscr G(t,x,v)/v ]\nonumber\\
&- |\nabla v|^4/v^2
+ \alpha(t) \partial_t v |\nabla v|^2/v^2 \nonumber\\
=&~ 2[\alpha(t)+1] \partial_t v |\nabla v|^2/v^2 
- 3 |\nabla v|^4/v^2 - \alpha(t) [\partial_t v/v]^2\nonumber\\
& -2 |\nabla v|^2 \mathscr G(t,x,v)/v^2
+\alpha (t) \partial_t v \mathscr G(t,x,v)/v^2.
\end{align}
 
Hence putting \eqref{eq-2.17}--\eqref{eq-2.18} together and rearranging terms it follows after 
some tedious but straightforward calculations that  
\begin{align}\label{eq-2.21}
2&(p-1) v \{ \langle \nabla ( |\nabla v|^2/v), \nabla v/v  \rangle
- \alpha(t) \langle\nabla ( \partial_t v/v), \nabla v/v \rangle \nonumber\\
&+\alpha(t) \langle \nabla ( \mathscr G(t,x,v)/v), \nabla v/v \rangle \}
+2 \langle \nabla v, \nabla |\nabla v|^2 \rangle/v
- 2\alpha(t) \langle \nabla v , \nabla \partial_t v \rangle/v\nonumber\\
&+ [ 2 |\nabla v|^2/v - \alpha(t) \partial_t v/v ] [ (p-1) \Delta_\phi v ]-|\nabla v|^4/v^2
+ \alpha(t) \partial_t v |\nabla v|^2/v^2\nonumber\\
=&~ 2p \langle \nabla F, \nabla v \rangle - [\partial_t v/v-|\nabla v|^2/v -\mathscr G(t,x,v)/v ]^2
+[1-\alpha (t)] [\partial_t v/v ]^2 \nonumber\\
&+ 2 \alpha (t) |\nabla v|^2 \mathscr G(t,x,v)/v^2 
+ [\mathscr G(t,x,v)/v]^2 -2 \alpha (t) \langle \nabla v ,\nabla \mathscr G(t,x,v) \rangle/v\nonumber\\
&+[\alpha (t)-2] \partial_t v \mathscr G(t,x,v)/v^2 \nonumber\\
=&~ 2p \langle \nabla F, \nabla v \rangle - [(p-1) \Delta_\phi v]^2 
+[1-\alpha (t)] [\partial_t v/v ]^2+ 2 \alpha (t) |\nabla v|^2 \mathscr G(t,x,v)/v^2 \nonumber\\
&+ [\mathscr G(t,x,v)/v]^2 - 2\alpha (t) \langle \nabla v ,\nabla \mathscr G(t,x,v) \rangle/v
+[\alpha (t)-2] \partial_t v \mathscr G(t,x,v)/v^2. 
\end{align}
Now substituting \eqref{eq-2.21} back into \eqref{eq-2.20} and simplifying terms results in  
\begin{align}\label{eq-2.26}
\mathscr L^p_v [F] =&~
2[\alpha(t) -1] h ( \nabla v, \nabla v)/v
- 2(p-1) |\nabla \nabla v|^2
-2(p-1) {\mathscr Ric}_\phi^m (\nabla v, \nabla v) \nonumber\\
&-2(p-1)\langle \nabla \phi , \nabla v \rangle^2/(m-n)
+ 2p \langle \nabla F, \nabla v \rangle- [(p-1) \Delta_\phi v]^2 \nonumber\\
&+2\alpha(t) (p-1) \langle h , \nabla \nabla v \rangle 
+ \alpha(t) (p-1) \langle 2 {\rm div} h - \nabla ({\rm Tr}_g h),\nabla v \rangle\nonumber\\
&+ \alpha(t) (p-1) \langle \nabla \partial_t \phi , \nabla v \rangle
- 2 \alpha(t) (p-1)h (\nabla \phi , \nabla v)+[1-\alpha (t)] [\partial_t v/v]^2\nonumber\\
&+2[1-\alpha(t)] \langle \nabla v ,\nabla \mathscr G(t,x,v) \rangle/v - \alpha (t)(p-1) \Delta_\phi [\mathscr G(t,x,v)]\nonumber\\
&+[\alpha(t)-1] |\nabla v|^2\mathscr G(t,x,v)/v^2
+ 2[\alpha(t)-1] \partial_t v \mathscr G(t,x,v)/v^2 \nonumber\\
&+[1-\alpha (t)] [\mathscr G(t,x,v)/v]^2
+ \alpha'(t) \mathscr G(t,x,v)/v -\alpha'(t) \partial_t v/v-\beta' (t).
\end{align}
Finally, making note of the relation 
\begin {align} \label{eq2.23}
[1-\alpha (t)] \left[\frac{\partial_t v}{v} \right]^2
&+ 2[\alpha(t)-1] \frac{\partial_t v}{v}  \frac{\mathscr G(t,x,v)}{v} 
+[1-\alpha (t)] \left[\frac{\mathscr G(t,x,v)}{v}\right]^2\nonumber\\
& =[1-\alpha (t)] \left[ \frac{\partial_t v}{v} - \frac{\mathscr G(t,x,v)}{v} \right]^2,
\end{align}
and substituting into \eqref{eq-2.26} gives the desired conclusion. 
\end{proof}

Having obtained the evolution identity for $F$, we now proceed on to deriving a corresponding inequality 
by applying the geometric bounds. As this is quite an involved task, for the sake of clarity we carry it out 
in a number of steps in the next few lemmas.

\begin{lemma}\label{lem-5.7}
Under the assumptions of Lemma $\ref{lemma-5.6-Equality}$ and subject to $\alpha(t) \ge 1$ we have 
\begin{align}\label{eq5.44}
\mathscr L^p_v [F] \le&
- \frac{1+m(p-1)}{m(p-1)}[(p-1)\Delta_\phi v]^2
+\frac{2}{v}[\alpha(t) -1] h ( \nabla v, \nabla v)\nonumber\\
&-2(p-1) {\mathscr Ric}_\phi^m (\nabla v, \nabla v) + 2p \langle \nabla v, \nabla F \rangle
+ (p-1) \alpha^2 (t) |h|^2  \nonumber\\
&+ \alpha(t) (p-1) \langle 2 {\rm div} h - \nabla ({\rm Tr}_g h),\nabla v \rangle
+ \alpha(t) (p-1) \langle \nabla \partial_t \phi , \nabla v \rangle\nonumber\\
&+\frac{2}{v}[1-\alpha(t)] \langle \nabla v, \nabla \mathscr G(t,x,v) \rangle
- \alpha (t)(p-1) \Delta_\phi [\mathscr G(t,x,v)] \nonumber\\
&- 2 \alpha(t) (p-1)h (\nabla \phi , \nabla v)
+[\alpha(t)-1]  \frac{|\nabla v|^2}{v} \frac{\mathscr G(t,x,v)}{v} \nonumber\\
&+ \alpha'(t) \frac{\mathscr G(t,x,v)}{v}
-\alpha'(t) \frac{\partial_t v}{v}-\beta' (t).
\end{align}
\end{lemma}

\begin{proof}
We start with \eqref{eq-5.27} in Lemma \ref{lemma-5.6-Equality}. 
As $[1-\alpha (t)] |\partial_t v/v - \mathscr G(t,x,v)/v|^2$ is non-positive 
($\alpha(t)\ge 1$) it can be omitted from the right-hand side thus giving an inequality. 
Next since by an application of the Cauchy-Schwarz inequality we have
\begin {align} \label{eq-5.47h}
\alpha(t)\langle h , \nabla \nabla v \rangle \le |\alpha(t)\langle h , \nabla \nabla v \rangle| 
\leq \frac{1}{2} |\nabla \nabla  v|^2 + \frac{1}{2} \alpha^2 (t) |h|^2, 
\end{align} 
and clearly $|\nabla \nabla v|^2+\langle \nabla \phi , \nabla v \rangle^2/(m-n) 
\ge (\Delta v)^2/n+\langle \nabla \phi, \nabla v \rangle^2/(m-n) 
\ge (\Delta_\phi v)^2/m$ we can write  
\begin{align}\label{eq-5.48}
\frac{(\Delta_\phi  v)^2}{m} - \alpha^2 (t) |h|^2 
& \le |\nabla \nabla  v|^2 + \frac{\langle \nabla \phi , \nabla v \rangle^2}{m-n} - \alpha^2 (t) |h|^2 \nonumber \\
& \le 2|\nabla \nabla v|^2 +2\frac{\langle \nabla \phi , \nabla v \rangle^2}{m-n} -2\alpha(t) \langle h , \nabla \nabla v \rangle. 
\end{align}

Substituting \eqref{eq2.23}, \eqref{eq-5.48} back into \eqref{eq-5.27} and taking into account the
condition $p>1$ gives the desired conclusion. 
\end{proof}

\begin{remark} \label{remark-on-b}
A close inspection of the above proof shows that when the metric is time independent 
(i.e., $\partial_t g=2h \equiv 0$) in the first term on the right-hand side of \eqref{eq5.44} 
the factor $[1+m(p-1)]/[m(p-1)]$ can be replaced with $[2+m(p-1)]/[m(p-1)]$ since the 
bound \eqref{eq-5.47h} is no longer needed.
\end{remark}

\begin{lemma} \label{khosh-akhlagh}
Under the assumptions of Lemma $\ref{lemma-5.6-Equality}$, 
with $b = [m(p-1)]/[1+ m(p-1)]$ and $\alpha(t) \ge 1$ we have 
\begin{align}\label{eq-5.51}
\mathscr L^p_v [F]
\le& -\frac{1}{b \alpha^2} F^2
 -\frac{2(\alpha -1)}{b \alpha^2} \frac{|\nabla v|^2}{v} F
 +\left[\mathscr G_v -\frac{2 \beta}{b \alpha^2} +\frac{\alpha'}{\alpha}\right]F
+ 2p \langle \nabla v, \nabla F \rangle \nonumber\\
&-\frac{(\alpha -1)^2}{b \alpha^2} \frac{|\nabla v|^4}{v^2}
+\frac{2}{v}(\alpha -1) h ( \nabla v, \nabla v)+ (p-1) \alpha^2 |h|^2 \nonumber\\
&+\bigg\{ (\alpha-1) \left[ \frac{\mathscr G}{v} - \mathscr G_v \right] 
- \alpha (p-1) v\mathscr G_{vv} - \frac{\alpha'}{\alpha} -\frac{2(\alpha -1)}{b \alpha^2} \beta \bigg\} 
\frac{|\nabla v|^2}{v} \nonumber \\
& -2(p-1) {\mathscr Ric}_\phi^m (\nabla v, \nabla v)
+ \alpha (p-1) \langle 2 {\rm div} h - \nabla ({\rm Tr}_g h),\nabla v \rangle\nonumber\\
&+2 \left[ (\alpha-1)\frac{|\mathscr G_x|}{v} + \alpha (p-1) |\mathscr G_{x v}| \right] |\nabla v|
+ \alpha (p-1) \langle \nabla \partial_t \phi , \nabla v \rangle\nonumber\\
&- 2 \alpha (p-1)h (\nabla \phi , \nabla v) - \alpha (p-1) \Delta_\phi \mathscr G^x
+ \beta \mathscr G_v -\frac{\beta}{\alpha} \left[ \frac{\beta}{b \alpha}
- \alpha'\right] -\beta'.
\end{align}
\end{lemma}

\begin{proof}
We start with some calculations relating to the nonlinearity $\mathscr G=\mathscr G(t,x,v)$,
abbreviating the arguments $(t, x, v)$ for convenience. Firstly, it is seen by calculating 
in local coordinates or directly that 
\begin {align} \label{eq6.44}
\nabla \mathscr G = \mathscr G_x +\nabla v \mathscr G_v, \qquad \mathscr G_x=(\mathscr G_{x_1}, \dots, \mathscr G_{x_n}). 
\end{align}  
Next, writing $\mathscr G^x: x \mapsto \mathscr G(t,x, v)$ we can differentiate \eqref{eq6.44} further to obtain  
\begin{align} \label{eq6.47}
\Delta \mathscr G &= \mathscr G_{x x} + 2\langle \nabla v,\mathscr G_{x v}\rangle 
+ \mathscr G_{vv} |\nabla v|^2+ \mathscr G_v \Delta v\nonumber\\
&= \Delta \mathscr G^x + 2\langle \nabla v ,\mathscr G_{x v} \rangle + \mathscr G_{vv} |\nabla v|^2+ \mathscr G_v \Delta v.
\end{align}   
For $\Delta_\phi \mathscr G$, by using the above calculations and substituting accordingly, it is seen that
\begin {align} \label{eq6.48}
\Delta_\phi \mathscr G &= \Delta \mathscr G - \langle \nabla \phi , \nabla \mathscr G \rangle 
= \Delta \mathscr G- \langle \nabla \phi, (\mathscr G_x + \mathscr G_v \nabla v) \rangle\nonumber\\
&= \Delta \mathscr G^x+ 2\langle \nabla v ,\mathscr G_{x v}\rangle + \mathscr G_{vv} |\nabla v|^2+ \mathscr G_v \Delta v
- \langle \nabla \phi, \mathscr G _x \rangle  - \mathscr G _v \langle \nabla \phi, \nabla v \rangle \nonumber\\
& = \Delta_ \phi \mathscr G^x + 2\langle \nabla v ,\mathscr G_{x v} \rangle + \mathscr G_{vv} |\nabla v|^2+ \mathscr G_v \Delta_\phi v.
\end{align}

We now use these calculations to simplify and rearrange the right-hand side of \eqref{eq5.44} in Lemma \ref{lem-5.7}. 
Towards this end, referring to the fourth line we have
\begin{align}\label{eq-5.56}
(2/v)[1-\alpha] \langle& \nabla v, \nabla \mathscr G \rangle
- \alpha (p-1) \Delta_\phi \mathscr G \nonumber\\
=&~ (2/v) [1-\alpha] \langle \nabla v, \mathscr G_x +\nabla v \mathscr G_v \rangle\nonumber\\
&- \alpha (p-1) [\Delta_ \phi \mathscr G^x + 2\langle \nabla v ,\mathscr G_{x v} \rangle + \mathscr G_{vv} |\nabla v|^2+ \mathscr G_v \Delta_\phi v]\nonumber\\
=&~ (2/v) \langle \nabla v, [1-\alpha] \mathscr G_x -\alpha (p-1) v\mathscr G_{xv} \rangle
 + (2/v) [1-\alpha] \mathscr G_v |\nabla v|^2 \nonumber\\
 & - \alpha (p-1) [\Delta_ \phi \mathscr G^x  + \mathscr G_{vv} |\nabla v|^2+ \mathscr G_v \Delta_\phi v].
\end{align}
Next, upon recalling \eqref{eq-8.1} and \eqref{eq5.21} we can write  
\begin{align} 
& [F - |\nabla v|^2/v +\beta]/\alpha= \mathscr G(t,x,v)/v-\partial_t v/v, \label{Eq-9.38.n} \\
& (p-1)\Delta_\phi v = [\partial_t v -|\nabla v|^2 - \mathscr G]/v
=- [F + (\alpha -1)|\nabla v|^2/v + \beta ]/\alpha. \label{eq-5.57} 
\end{align}
Substituting \eqref{eq-5.56}, \eqref{Eq-9.38.n} and \eqref{eq-5.57} 
back into the right-hand side of \eqref{eq5.44}, we can re-write the inequality as 
\begin{align} \label{sardard}
\mathscr L^p_v [F] \le& 
-[1/(b \alpha^2)] [F + (\alpha -1) |\nabla v|^2/v + \beta ]^2 
\nonumber \\
&+ 2p \langle \nabla v, \nabla F \rangle 
-2(p-1) {\mathscr Ric}_\phi^m (\nabla v, \nabla v) \nonumber\\
&+ (|\nabla v|^2/v) \{ (\alpha -1) [\mathscr G -2v \mathscr G_v]/v
- \alpha (p-1) v\mathscr G_{vv} \} \nonumber\\
&+(2/v)(\alpha -1) h ( \nabla v, \nabla v)
+ \alpha (p-1) \langle 2 {\rm div} h - \nabla ({\rm Tr}_g h),\nabla v \rangle\nonumber\\
&+ \alpha (p-1) \langle \nabla \partial_t \phi , \nabla v \rangle 
- 2 \alpha (p-1)h (\nabla \phi , \nabla v)
- \alpha (p-1) \Delta_\phi \mathscr G^x \nonumber\\
&+(2/v) \langle \nabla v, (1-\alpha) \mathscr G_x -\alpha (p-1)v \mathscr G_{xv} \rangle
 + (p-1) \alpha^2 |h|^2-\beta' \nonumber\\
&+ \mathscr G_v [F + (\alpha -1)|\nabla v|^2/v + \beta ] 
+\alpha'[F- |\nabla v|^2/v+\beta]/\alpha.
\end{align}

Referring to the right-hand side in the inequality above, 
the sum of the expression on the first line and the two terms 
on the last line can be expanded and rewritten as  
\begin{align}
-[1/(b \alpha^2)] &[F + (\alpha -1) |\nabla v|^2/v + \beta]^2 \nonumber\\
+ \mathscr G_v &[F + (\alpha -1)|\nabla v|^2/v + \beta]+ \alpha' [F- |\nabla v|^2/v+\beta]/\alpha \nonumber\\
=&-[1/(b \alpha^2)]\{F^2 +(\alpha -1)^2 |\nabla v|^4/v^2 + \beta^2
+ 2(\alpha -1) (|\nabla v|^2 /v) F\nonumber\\
&+ 2 \beta F + 2(\alpha -1) \beta |\nabla v|^2/v \} + \mathscr G_v F
+ (\alpha -1)  \mathscr G_v|\nabla v|^2/v \nonumber\\
&+ \beta \mathscr G_v + \alpha'/\alpha F
- \alpha' |\nabla v|^2/(\alpha v)+ \alpha' \beta/\alpha\nonumber\\
=&-[1/(b \alpha^2)] F^2 -2(\alpha -1)/(b \alpha^2)(|\nabla v|^2/v)F
+\{\alpha'/\alpha - 2 \beta/(b \alpha^2) + \mathscr G_v\}F\nonumber\\
&-(\alpha -1)^2 /(b \alpha^2) (|\nabla v|^4/v^2)
-(\beta/\alpha) \{\beta/(b \alpha)- \alpha'\}+ \beta \mathscr G_v\nonumber\\
&- \{2(\alpha -1)\beta/(b \alpha^2)
+\alpha'/\alpha-(\alpha -1) \mathscr G_v \}(|\nabla v|^2/v).
\end{align}
Substituting the above back in \eqref{sardard} gives 
\begin{align} \label{sarkhoob}
\mathscr L^p_v [F] \le& 
-[1/(b \alpha^2)]F^2 
- 2(\alpha -1) /(b \alpha^2) (|\nabla v|^2/v)F
+\{\alpha'/\alpha - 2 \beta/(b \alpha^2) + \mathscr G_v \}F\nonumber\\
&+ 2p \langle \nabla v, \nabla F \rangle- (\alpha -1)^2 /(b \alpha^2) ( |\nabla v|^4/v^2)
+(2/v)(\alpha -1) h ( \nabla v, \nabla v)\nonumber\\
&-2(p-1) {\mathscr Ric}_\phi^m (\nabla v, \nabla v) 
+ \{ (\alpha -1) [\mathscr G-2v \mathscr G_v]/v - \alpha (p-1) v\mathscr G_{vv}\nonumber\\
& + (\alpha -1) \mathscr G_v-2(\alpha -1)\beta/(b \alpha^2)
- \alpha'/\alpha \} (|\nabla v|^2/v)\nonumber\\
&+ \alpha (p-1) \langle 2 {\rm div} h - \nabla ({\rm Tr}_g h),\nabla v \rangle
+ \alpha (p-1) \langle \nabla \partial_t \phi , \nabla v \rangle\nonumber\\
& - 2 \alpha (p-1)h (\nabla \phi , \nabla v)+ (p-1) \alpha^2 |h|^2
- \alpha (p-1) \Delta_\phi \mathscr G^x \nonumber\\
&+(2/v) \langle \nabla v, (1-\alpha) \mathscr G_x -\alpha (p-1)v \mathscr G_{xv} \rangle-\beta' \nonumber\\
&-(\beta/\alpha) \{ \beta/(b \alpha)- \alpha'\} + \beta \mathscr G_v.
\end{align}

Next, referring to the first term in the penultimate line on the right-hand side of 
\eqref{sarkhoob}, by using the Cauchy-Schwarz inequality, 
we have the following upper bound 
\begin{align} \label{sarkhoob-pb1}
(2/v)\langle \nabla v, (1-\alpha) \mathscr G_x -\alpha (p-1) v \mathscr G_{x v} \rangle 
&\le (2|\nabla v|/|v|)[ |1-\alpha| |\mathscr G_x| +  \alpha |p-1| |v| | \mathscr G_{x v}|] \nonumber\\
&\le 2 |\nabla v| [(\alpha-1) |\mathscr G_x|/v + \alpha (p-1) |\mathscr G_{x v}|]
\end{align} 
as $p>1$ and $\alpha \ge 1$. Now making use of the upper bound 
\eqref{sarkhoob-pb1} in \eqref{sarkhoob} gives at once the desired conclusion. 
\end{proof}

\begin{lemma}\label{lem-2.3}
Under the assumptions of Lemma $\ref{lemma-5.6-Equality}$, 
with $b = [m(p-1)]/[1+ m(p-1)]$, $\alpha(t) \ge 1$ and the bounds in Theorem $\ref {thm5.1}$, 
in $Q_{2R, T}$ we have  
\begin{align}\label{eq-5.53}
\mathscr L^p_v [F]
\le& -\frac{F^2}{b \alpha^2}  
-\frac{2(\alpha -1)}{b \alpha^2} \frac{|\nabla v|^2}{v} F
+\left[\frac{\alpha'}{\alpha} -\frac{2 \beta}{b \alpha^2} + \mathscr G_v\right]F
+ 2p \langle \nabla v, \nabla F \rangle \\
&-\frac{(\alpha -1)^2}{b \alpha^2} \frac{|\nabla v|^4}{v^2}
+\bigg\{2(p-1)v [(m-1) k+k_2]+2 (\alpha -1) \overline k_1\nonumber\\
&+(\alpha-1)\left[ \frac{\mathscr G}{v} - \mathscr G_v \right]
- \alpha (p-1) v\mathscr G_{vv} -\frac{2(\alpha -1)}{b \alpha^2} \beta
- \frac{\alpha'}{\alpha}\bigg\}\frac{|\nabla v|^2}{v}\nonumber\\
&+\bigg\{ 2 (\alpha-1)\frac{|\mathscr G_x|}{v} + 2\alpha (p-1)  |\mathscr G_{x v}|
+ \alpha (p-1)  \ell_2+ 2 \alpha (p-1) \underline k_1 \ell_1 \bigg\}|\nabla v|\nonumber\\
&+ \alpha^2 (p-1)n[(\underline k_1 + \overline k_1)^2+2 k_2] - \alpha (p-1) \Delta_\phi \mathscr G^x
+ \beta \mathscr G_v -\frac{\beta}{\alpha} \left[ \frac{\beta}{b \alpha}
- \alpha'\right] -\beta'.\nonumber
\end{align}

\end{lemma}

\begin{proof}
We start with the inequality established in Lemma \ref{khosh-akhlagh} and apply the bounds 
in Theorem \ref{thm5.1}. Indeed from \eqref {vahideh-A-2} we firstly note  
$\alpha^2 |h|^2  \le \alpha^2 (\underline k_1 + \overline k_1)^2 |g|^2 
= n \alpha^2 (\underline k_1 + \overline k_1)^2$, and likewise upon recalling the relation 
${\partial_t} g= 2h$, that 
\begin{align}\label{EQ-eq-3.45}
-2\underline k_1 |\nabla v|^2 \le 2 h ( \nabla v, \nabla v) \le 2 \overline k_1 |\nabla v|^2.
\end{align} 
Next using $|g^{ij} (2 \nabla_i h_{j \ell} - \nabla_\ell h_{ij})| 
\le 3|g| |\nabla h|$ along with $2|\nabla h| = |\nabla \partial_t g| \le 2k_2$ 
we can write after an application of the Cauchy-Schwarz inequality that
\begin {align*} 
\langle  {\rm div} h - \frac{1}{2} \nabla ({\rm Tr}_g h),\nabla v \rangle
&\le \left|{\rm div} h - \frac{1}{2} \nabla ({\rm Tr}_g h)\right| |\nabla v| 
= \left| g^{ij} \nabla_i h_{j \ell} - \frac{1}{2} g^{ij} \nabla_\ell h_{ij} \right| |\nabla v|.
\end{align*}
Therefore, by an application of Young's inequality it is seen that 
\begin{align}\label{EQ-eq-3.46}
\alpha \langle  2 {\rm div} h - \nabla ({\rm Tr}_g h),\nabla v \rangle
&\le \alpha \left| g^{ij} (2 \nabla_i h_{j \ell} - \nabla_\ell h_{ij}) \right| |\nabla v| \nonumber\\
&\le 3 \alpha \sqrt n k_2 |\nabla v| \nonumber\\
&\le 2 \alpha^2 n k_2 + 2 k_2 |\nabla v|^2.
\end{align} 

Finally, by virtue of the bounds $|\nabla \phi| \le \ell_1$ and $|\nabla \partial_t \phi| \le \ell_2$, 
upon an application of the Cauchy-Schwarz inequality we have
\begin{align}\label{EQ-eq-3.47}
\alpha \langle \nabla \partial_t \phi , \nabla v \rangle
\le \alpha |\nabla \partial_t \phi| |\nabla v|
\le \alpha \ell_2 |\nabla v|,
\end{align}
and in a similar way 
\begin{align}\label{EQ-eq-3.48}
- \alpha \underline k_1 \ell_1 |\nabla v| \le -\alpha \underline k_1 |\nabla \phi| |\nabla v| 
\le \alpha h (\nabla \phi , \nabla v) \le \alpha \overline k_1 |\nabla \phi| |\nabla v|
\le \alpha \overline k_1 \ell_1 |\nabla v|.
\end{align}
Substituting the above in the inequality \eqref{eq-5.51} established in Lemma \ref{khosh-akhlagh} 
and making note of the conditions $p>1$ with $\alpha \ge 1$ leads to the desired conclusion. 
\end{proof}

\section{Proof of the local estimate in Theorem \ref{thm5.1}}
\label{sec4}

The proof of the Li-Yau estimate is based on localisation and the results from the 
previous section. A key ingredient in localisation is the construction and use of suitable cut-off 
functions. The following standard lemma gives the main properties for a profile function that is needed later.

\begin{lemma} \label{psi lemma} 
There exists a function $\bar{\eta}:[0,\infty) \to \mathbb{R}$ satisfying the following properties:
\begin{enumerate}[label=$(\roman*)$]
\item $\bar\eta$ is of class $\mathscr{C}^2 [0, \infty)$. 
\item $0 \le \bar\eta(s) \le 1$ for $0 \le s < \infty$ and $\bar{\eta} \equiv 1$ on $[0,1]$ and $\bar{\eta} \equiv 0$ on $[2, \infty)$.
\item $\bar\eta' \le 0$ and so $\bar\eta$ is non-increasing and for suitable constants $c_1, c_2>0$ we have
\begin{align} \label{9.30}
- c_1  \bar\eta^{1/2} \le \bar\eta' \le 0, \qquad and \qquad  \bar{\eta}^{''} \ge -c_2. 
\end{align}
\end{enumerate}
\end{lemma}

Now we pick a reference point $x_0 \in M$, fix $R, T>0$ and $0<T_1 \le T$ and then with $\varrho (x,t)$ denoting the 
geodesic radial variable with respect to $x_0$ at time $t$, set
\begin{align}\label{9.31}
\eta(x,t)=\bar{\eta} (\varrho (x, t)/R), \qquad x \in M, \quad 0 \le t \le T. 
\end{align}
It is evident that the resulting function $\eta$ satisfies $\eta \equiv 1$ for when $0 \le \varrho(x, t) \le R$ and $\eta \equiv 0$ for 
when $\varrho(x,t) \ge 2R$. Additionally from 
\eqref{9.31} we have  
$\nabla \eta = (\bar\eta'/R) \nabla \varrho$ and $\Delta \eta = \bar\eta'' |\nabla \varrho|^2/R^2+\bar\eta' \Delta \varrho/R$ 
and so 
$\Delta_\phi \eta = \Delta \eta - \langle \nabla \phi, \nabla \eta \rangle = \bar\eta'' |\nabla \varrho|^2/R^2+\bar\eta' \Delta_\phi \varrho/R$. 
In particular $\nabla \eta$, $\Delta_\phi \eta$ vanish outside the space-time set $R \le \varrho(x,t) \le 2R$. Finally, by putting the above 
fragments together, it is easily see that 
\begin{align}
\mathscr L^p_v [\eta] &= [\partial _t - (p-1) v \Delta_ \phi] \eta = \bar \eta' \partial_t \varrho/R 
- (p-1) v [\bar\eta'' |\nabla \varrho|^2/R^2+\bar\eta' \Delta_\phi \varrho/R] \nonumber \\
&= (\bar\eta'/R) [\partial _t - (p-1) v \Delta_ \phi] \varrho - (p-1) v (\bar\eta''/R^2)  |\nabla \varrho|^2 \nonumber \\
&= (\bar\eta'/R) \mathscr L^p_v [\varrho] - (p-1) v (\bar\eta''/R^2)  |\nabla \varrho|^2.
\end{align}

\qquad \\
Let us now move on to the proof of Theorem \ref{thm5.1}. 
Towards this end consider the spatially localised function $G= t\eta F$ 
where $F$ is as in \eqref{eq5.21}. Let $(x_1,t_1)$ denote the point where $G$ attains its maximum over the compact 
cylinder $\{\varrho(x,t) \le 2R, 0 \le t \le T_1\}$. We assume $G (x_1, t_1)>0$ as otherwise the estimate immediately follows 
from $G= t\eta F \le 0$. So in particular $t_1>0$ and from $(ii)$ in Lemma \ref{psi lemma} it must be that $\varrho (x_1, t_1) < 2R$. 
Therefore by basic considerations at the maximum point $(x_1, t_1)$ we have 
\begin{align} \label{9.32}
\partial_t G \ge 0, \qquad
\nabla G =0, \qquad
\Delta_\phi G \le 0.
\end{align}
Using the identity 
$\langle \nabla \eta, \nabla (\eta F) \rangle = \eta \langle \nabla \eta, \nabla F \rangle + F |\nabla \eta|^2$ we can write  
\begin{align} 
\mathscr L^p_v [G] 
=&~\mathscr L^p_v[tF \eta] = t F \mathscr L^p_v [\eta] - 2 t(p-1)v\langle \nabla \eta, \nabla F \rangle
+ t \eta \mathscr L^p_v [F] + \eta F\\
=&~  t F \mathscr L^p_v [\eta] - 2 t(p-1)\langle \nabla \eta, \nabla (\eta F) \rangle/\eta + 2 t (p-1)v F |\nabla \eta|^2/\eta
+ t \eta \mathscr L^p_v [F] + \eta F.\nonumber
\end{align}
Therefore, utilising $\nabla G=0$ and $\mathscr L^p_v[G]= [\partial_t - (p-1)v \Delta_\phi]G\ge 0$ 
from \eqref{9.32}, we have at the maximum point 
$(x_1, t_1)$, the inequality  
\begin{align}  \label{9.34}
0\le \mathscr L^p_v [G] 
=&~  t_1 F \mathscr L^p_v [\eta]+ 2 t_1 (p-1)v F |\nabla \eta|^2/\eta
+ t_1 \eta \mathscr L^p_v [F] + \eta F.
\end{align}
Substituting from \eqref{eq-5.53} into \eqref{9.34} then gives   
\begin{align}\label{eq-2.38}
0 \le&~ \mathscr L^p_v[G] =
t_1 \eta \mathscr L^p_v [F] + 2  t_1 (p-1) v\frac{|\nabla \eta |^2}{\eta} F 
+ t_1 F \mathscr L^p_v [\eta] + \eta F \nonumber\\
\le& -\frac{t_1 \eta}{b \alpha^2} F^2 
- \frac{2(\alpha -1)}{b \alpha^2} \frac{|\nabla v|^2}{v} t_1 \eta F
+\left[\frac{\alpha'}{\alpha} -\frac{2 \beta}{b \alpha^2}+ \mathscr G_v\right]t_1 \eta F
+ 2p t_1 \eta \langle \nabla v, \nabla F \rangle \nonumber\\
&- \frac{(\alpha -1)^2}{b \alpha^2}t_1 \eta \frac{|\nabla v|^4}{v^2} 
+ \bigg\{2(p-1)v [(m-1) k+k_2]+2 (\alpha -1) \overline k_1\nonumber\\
&+(\alpha-1)\left[ \frac{\mathscr G}{v} - \mathscr G_v \right]
- \alpha (p-1) v\mathscr G_{vv} -\frac{2(\alpha -1)}{b \alpha^2} \beta
- \frac{\alpha'}{\alpha}\bigg\} t_1 \eta \frac{|\nabla v|^2}{v}\nonumber\\
&+ \bigg\{2 (\alpha-1) \frac{|\mathscr G_x|}{v}+ 2\alpha (p-1)  |\mathscr G_{x v}| 
+ \alpha (p-1)  \ell_2+ 2 \alpha (p-1) \underline k_1 \ell_1 \bigg\} t_1 \eta |\nabla v|\nonumber\\
&+ \alpha^2 (p-1)t_1 \eta n[(\underline k_1 + \overline k_1)^2+2k_2]
+ \beta t_1 \eta \mathscr G_v - \alpha (p-1) t_1 \eta \Delta_\phi \mathscr G^x \nonumber\\
&- t_1 \eta \frac{\beta}{\alpha}\left[ \frac{\beta}{b \alpha}- \alpha'\right] 
 - t_1 \eta \beta' + 2t_1(p-1) v \frac{|\nabla \eta|^2}{\eta} F\nonumber\\
 &+ t_1 F \partial_t \eta- t_1F(p-1)v \Delta_\phi \eta + \eta F.
\end{align}

Multiplying through by $\eta$, using $\langle \nabla v, \nabla F \rangle \le |\nabla v| |\nabla F| 
= |\nabla v| (|\nabla \eta|/\eta) F$ [the last identity being a result of $\eta \nabla F = -F\nabla \eta $ that holds 
at the maximum point $(x_1, t_1)$] and re-writing the expression on the right-hand side in terms of $G$, 
 \begin{align}\label{eq2.46}
0 \le& -\frac{G^2}{b \alpha^2t_1} 
- \frac{2(\alpha -1)}{b \alpha^2}  \frac{|\nabla v|^2}{v} \eta G
+ 2p |\nabla v| \frac{|\nabla \eta|}{\eta} \eta G
- \frac{(\alpha -1)^2}{b \alpha^2} t_1\eta^2 \frac{|\nabla v|^4}{v^2} \nonumber\\
&+\left\{\frac{\partial_t \eta}{\eta}-(p-1)v \frac{\Delta_\phi \eta}{\eta} + 2(p-1) v \frac{|\nabla \eta|^2}{\eta^2}
+\frac{\alpha'}{\alpha} -\frac{2 \beta}{b \alpha^2}+ \mathscr G_v +\frac{1}{t_1} \right\} \eta G\nonumber\\
&+ \bigg\{2(p-1)v [(m-1) k+k_2]+2 (\alpha -1) \overline k_1
+(\alpha -1)\left[\frac{\mathscr G}{v} -\mathscr G_v \right]\nonumber\\
&- \alpha(p-1) v\mathscr G_{vv} -\frac{2(\alpha -1)}{b \alpha^2} \beta
 - \frac{\alpha'}{\alpha}\bigg\} t_1 \eta^2 \frac{|\nabla v|^2}{v}
+ \bigg\{2 (\alpha-1) \frac{|\mathscr G_x|}{v}\nonumber\\
&+ 2\alpha (p-1)  |\mathscr G_{x v}| + \alpha (p-1)  \ell_2
+ 2 \alpha (p-1) \underline k_1 \ell_1 \bigg\} t_1 \eta^2 |\nabla v|\nonumber\\
&-\bigg\{ \beta' + \alpha (p-1) \Delta_\phi \mathscr G^x - \beta \mathscr G_v
+ \frac{\beta}{\alpha}\left[ \frac{\beta}{b \alpha} -\alpha' \right] \nonumber\\
&- \alpha^2 (p-1)n[(\underline k_1 + \overline k_1)^2+2k_2] \bigg\} t_1 \eta^2. 
\end{align}

Let us now proceed by bounding from below the terms on the right-hand side of \eqref{eq2.46}. 
To this end we start with the expression on the second line involving $\mathscr L^p_v [\eta]/\eta$. To 
handle this term we consider the space and time derivatives separately and then combine the two at the end. 
\begin{itemize}
\item The term $(p-1)v \Delta_\phi \eta$: We assert that $\Delta_\phi \eta$ can be bounded from below as 
\begin{align} \label{eq-6.14}
\Delta_\phi \eta \ge -\frac{1}{R^2}[c_2 +(m-1) c_1(1+R \sqrt{k})], 
\end{align}
where $c_1, c_2>0$ are as in \eqref{9.30}. 
In fact, in view of ${\mathscr Ric}_\phi^m (g) \ge -(m-1) k g$, it follows from the weighted 
Laplacian comparison theorem that 
\begin{align}
\Delta_\phi \varrho \le (m-1) \sqrt {k}\coth (\sqrt{k} \varrho).
\end{align}
Using the $\phi$-Laplacian relation
$\Delta_\phi \eta = \bar\eta'' |\nabla \varrho|^2/R^2+\bar\eta' \Delta_\phi \varrho/R$, 
the bounds in Lemma \ref{psi lemma} and the fact that $\Delta_\phi \eta \equiv 0$ outside $R \le \varrho \le 2R$, we can write
\begin{align}\label{9.38}
\Delta_\phi \eta 
&\ge \frac{\bar\eta'' }{R^2}+ \frac {(m-1)}{R} \bar\eta' \sqrt{k} \coth (\sqrt{k} \varrho) 
\ge \frac{\bar\eta'' }{R^2} + \frac {(m-1)}{R} \bar\eta' \sqrt{k} \coth (\sqrt{k} R) \nonumber \\
&\ge \frac{ \bar\eta'' }{R^2} + (m-1) \frac{1+\sqrt{k}R}{R} \bar\eta' 
\ge -\frac{c_2}{R^2}-(m-1)\frac{c_1}{R} \left(\frac{1}{R}+\sqrt{k}\right).
\end{align} 
We have used the monotonicity of $s \mapsto \coth s$ and $s \coth s \le (1 + s)$ for $s>0$.

\item The term $\partial_t \eta$: Upon recalling \eqref{9.31}, it is evident that
$ \partial_t \eta = \bar\eta' \partial_t \varrho/R$.
Now utilising Lemma \ref{psi lemma} [the first inequality in \eqref{9.30}] and the geometric 
bound $\partial_t \varrho \ge - \underline k_1 \varrho$, this gives 
\begin{align} \label{9.41}
 \partial_t \eta 
=  \bar\eta' \partial_t \varrho/R 
\le -\varrho \underline k_1 \bar\eta'/R  \le c_1 \varrho \underline k_1 \sqrt{\bar\eta} /R  \le  2c_1 \underline k_1. 
\end{align}
Here in the last inequality we have used $0\le \varrho \le 2R$.
The justification for the lower bound $\partial_t \varrho(x_1,t_1) \ge - 2\underline k_1 R$ is as follows. 
Fix $(x,t)$ so that $d(x, x_0, t)<R$. Let $X(x_0,x)$ be the set 
of minimal geodesics $\lambda=\lambda(s): [0, 1] \to {\mathscr M}$ with respect to $g(t)$ connecting 
$x_0 =\lambda(0)$ to $x=\lambda(1)$ and let $\Omega(x_0,x)$ be the set of ${\mathscr C}^1$ 
curves connecting $x_0$ to $x$. Using $\partial_t g \ge -2 \underline k_1 g$ 
in $Q_{R,T}$ and Lemma {\rm B.40} p.~531 in \cite{Chow},
\begin{align} \label{rt-bound}
\partial_t \varrho (x,t) &= \frac{\partial}{\partial t} d (x, x_0; t) 
= \frac{\partial}{\partial t} \left\{ \inf_{\omega \in \Omega(x_0, x)} \int_0^1 |\omega'(s)|_{g(t)} \, ds \right\} \\
&= \frac{\partial}{\partial t} \left\{ \inf_{\omega \in \Omega(x_0, x)} \int_0^1 \sqrt{[g(t)] (\omega'(s), \omega'(s))} \, ds \right\} \nonumber \\
&= \inf_{\lambda \in X(x_0,x)} \int_0^1 \frac{[\partial_t g(t)](\lambda'(s), \lambda'(s))}{2  \sqrt{[g(t)] (\lambda'(s), \lambda'(s))}} \, ds \nonumber\\
&=  \inf_{\lambda \in X(x_0,x)} \int_0^1 \frac{[\partial_t g(t)](\lambda'(s), \lambda'(s))}{2|\lambda'(s)|_{g(t)}} \, ds 
\ge - \underline k_1 \varrho(x, t) \ge - 2\underline k_1 R. \nonumber 
\end{align}
\end{itemize}

In conclusion by putting together \eqref{eq-6.14} and \eqref{9.41} we arrive at the lower bound
\begin{align} \label{Lpv-Ali-Added}
\mathscr L^p_v [\eta] 
&\le 2c_1 \underline k_1 + (p-1) v \frac{1}{R^2}[c_2 +(m-1) c_1(1+R \sqrt{k})] . 
\end{align}

Moving to the third term in the second line on the right in \eqref{eq2.46}, again by utilising \eqref{9.30}
in Lemma \ref{psi lemma} and \eqref{9.31}, we have
\begin{align} \label{Eq5.35}
\frac{|\nabla \eta |^2}{\eta} = \frac{\bar\eta'^2}{\bar\eta} \frac{|\nabla \varrho|^2}{R^2} 
= \left( \frac{\bar\eta'}{\sqrt {\bar\eta}}  \frac{|\nabla \varrho|}{R} \right)^2 \le  \frac{c_1 ^2}{R^2}.  
\end{align}

Therefore, substituting the bounds \eqref{Lpv-Ali-Added} and \eqref{Eq5.35} back into 
the inequality \eqref{eq2.46} and rearranging terms gives 
\begin{align}\label{eq-6.15}
0 \le& -\frac{G^2}{b \alpha^2t_1} 
- \frac{2(\alpha -1)}{b \alpha^2}  \frac{|\nabla v|^2}{v} \eta G
+ 2p \sqrt v \frac{|\nabla \eta|}{\sqrt \eta} \frac{|\nabla v|}{\sqrt v} \sqrt \eta G
- \frac{(\alpha -1)^2}{b \alpha^2} t_1\eta^2 \frac{|\nabla v|^4}{v^2} \\
&+\left\{(p-1)v \left[ \frac{1}{R^2}[c_2 +(m-1) c_1(1+R \sqrt{k})]
+ 2 \frac{c_1^2}{R^2}\right]+ 2c_1 \underline k_1\right\}G \nonumber\\
&+\left[\frac{\alpha'}{\alpha} -\frac{2 \beta}{b \alpha^2} + \mathscr G_v + \frac{1}{t_1} \right] \eta G
+ \bigg\{2(p-1)v [(m-1) k+k_2]+2 (\alpha -1) \overline k_1\nonumber\\
&+(\alpha -1)\left[\frac{\mathscr G}{v} -\mathscr G_v \right]
- \alpha(p-1) v\mathscr G_{vv} -\frac{2(\alpha -1)}{b \alpha^2} \beta
- \frac{\alpha'}{\alpha}\bigg\} t_1 \eta^2 \frac{|\nabla v|^2}{v}\nonumber\\
&+ \bigg\{2 (\alpha-1) \frac{|\mathscr G_x|}{v}
+ 2\alpha (p-1)  |\mathscr G_{x v}| + \alpha (p-1)  \ell_2
+ 2 \alpha (p-1) \underline k_1 \ell_1 \bigg\} t_1 \eta^2 |\nabla v|\nonumber\\
&-\bigg\{ \beta' + \alpha (p-1) \Delta_\phi \mathscr G^x - \beta \mathscr G_v
+ \frac{\beta}{\alpha}\left[ \frac{\beta}{b \alpha} -\alpha' \right]
- \alpha^2 (p-1)n[(\underline k_1 + \overline k_1)^2+2k_2] \bigg\} t_1 \eta^2. \nonumber
\end{align}  

In order to bound the algebraic sum of the second and third terms on the right-hand side of \eqref{eq-6.15} 
we make use of the inequality $-\mathsf ax^2 +\mathsf bx \le \mathsf b^2 /(4\mathsf a)$ with $x$ being $|\nabla v|/\sqrt v$ and $\mathsf a > 0$, 
and taking advantage of \eqref{9.30} in Lemma \ref{psi lemma}, to write 
\begin{align}\label{eq5.68}
-\frac{2(\alpha -1)}{b \alpha^2}  \frac{|\nabla v|^2}{v} \eta G 
+ 2p \sqrt v \frac{|\nabla \eta|}{\sqrt \eta} \frac{|\nabla v|}{\sqrt v} \sqrt \eta G 
\le \frac{b \alpha^2 p^2 v}{2(\alpha -1)} \frac{c_1 ^2}{R^2} G.
\end{align}

Next, for the penultimate line on the right-hand side of \eqref{eq-6.15}, 
by an application of Young’s inequality, we can write
\begin{align}\label{eq5.72}
\bigg\{ & (\alpha-1) \frac{|\mathscr G_x|}{v}+ \alpha (p-1)  |\mathscr G_{x v}|+ 
\alpha (p-1)  \ell_2/2+ \alpha (p-1) \underline k_1 \ell_1 \bigg\}|\nabla v| \\
\le&~ \frac{3v^{2/3}}{4\varepsilon^{1/3}} \left\{ (\alpha-1) \frac{|\mathscr G_x|}{v}+ \alpha (p-1)  |\mathscr G_{x v}| 
+\alpha (p-1)  \ell_2/2 + \alpha (p-1) \underline k_1 \ell_1 \right\}^{4/3} + \varepsilon\frac{|\nabla v|^4}{4v^2}.\nonumber
\end{align}
Now substituting \eqref{eq5.68}--\eqref{eq5.72} back into \eqref{eq-6.15} and rearranging terms results in

\begin{align} \label{eq-2.44}
0 \le& -\frac{G^2}{b \alpha^2 t_1} 
+\frac{b \alpha^2 p^2 v}{2(\alpha -1)} \frac{c_1 ^2}{R^2} G
-\left[\frac{(\alpha -1)^2}{b \alpha^2}-\frac{\varepsilon}{2} \right] t_1\eta^2 \frac{|\nabla v|^4}{v^2} 
+\left[\frac{\alpha'}{\alpha} -\frac{2 \beta}{b \alpha^2}
+ \mathscr G_v +\frac{1}{t_1}\right]\eta G \nonumber\\
&+\left\{ (p-1)v \left[ \frac{1}{R^2}[c_2 +(m-1) c_1(1+R \sqrt{k})]
+ 2 \frac{c_1^2}{R^2}\right]+ 2c_1 \underline k_1 \right\} G\nonumber\\
&+\bigg\{2(p-1)v [(m-1) k+k_2]+ 2(\alpha -1) \overline k_1 
-\frac{2(\alpha -1)}{b \alpha^2} \beta - \frac{\alpha'}{\alpha} \nonumber\\
&+(\alpha -1)\left[\frac{\mathscr G}{v} -\mathscr G_v \right]
- \alpha(p-1) v\mathscr G_{vv}  \bigg\} t_1 \eta^2  \frac{|\nabla v|^2}{v} 
+ \bigg\{ - \beta' - \alpha (p-1) \Delta_\phi \mathscr G^x\nonumber\\
&+ \frac{3v^{2/3}}{2\varepsilon^{1/3}}\left[ (\alpha-1) \frac{|\mathscr G_x|}{v}+ \alpha (p-1)  |\mathscr G_{x v}| 
+\alpha (p-1)  \ell_2/2 + \alpha (p-1) \underline k_1 \ell_1 \right]^{4/3} \nonumber\\
& +\beta \mathscr G_v- \frac{\beta}{\alpha} \left[ \frac{\beta}{b \alpha} -\alpha' \right] 
+ \alpha^2 (p-1)n[(\underline k_1 + \overline k_1)^2+2k_2] \bigg\} t_1 \eta ^2.
\end{align}  

Using the inequality $-\mathsf ax^2 +\mathsf bx \le \mathsf b^2 /(4\mathsf a)$ 
with $x$ being $|\nabla v|^2/v$ and $\mathsf a > 0$ (with the improvement 
$-\mathsf ax^2 + \mathsf bx \le 0$ for $\mathsf b \le 0$) and noting $0<\varepsilon<2(\alpha-1)^2/(b \alpha^2)$, 
we have 
\begin{align}\label{eq-2.55}
-\frac{2(\alpha -1)^2- \varepsilon b\alpha^2}{2b \alpha^2}&
 t_1\eta^2 \frac{|\nabla v|^4}{v^2} + \bigg\{ 2(p-1)v [(m-1) k+k_2]
+ 2(\alpha -1) \overline k_1 \nonumber\\
& - \frac{2(\alpha -1)}{b \alpha^2} \beta
- \frac{\alpha'}{\alpha} + (\alpha -1)\left[\frac{\mathscr G}{v} -\mathscr G_v \right]
- \alpha(p-1) v\mathscr G_{vv}  \bigg\} t_1 \eta ^2 \frac{|\nabla v|^2}{v} \nonumber\\
\le&~ \frac{ b \alpha^2 t_1 \eta^2}{ 4 (\alpha -1)^2-2\varepsilon b \alpha^2} 
\bigg[ 2(p-1)v [(m-1) k+k_2] +2(\alpha -1) \overline k_1 \nonumber\\
& - \frac{2(\alpha -1)}{b \alpha^2} \beta- \frac{\alpha'}{\alpha}+ (\alpha -1)\left[\frac{\mathscr G}{v} -\mathscr G_v \right]
- \alpha(p-1) v\mathscr G_{vv} \bigg]_+^2.
\end{align}

Substituting \eqref{eq-2.55} back into \eqref{eq-2.44}, multiplying through by $t_1$ and rearranging terms, it follows that
\begin{align} \label{eq-2.49-pre}
0 \le& - \frac{1}{b \alpha^2} G^2
+ \bigg\{\frac{b \alpha^2 p^2 v}{2(\alpha -1)} \frac{c_1 ^2}{R^2}
+\left[\mathscr G_v +\frac{\alpha'}{\alpha} -\frac{2 \beta}{b \alpha^2} +\frac{1}{t_1} \right] \eta+ 2c_1 \underline k_1\nonumber\\
&+(p-1)v \left[\frac{1}{R^2}[c_2 +(m-1) c_1(1+R \sqrt{k})] 
+ 2 \frac{c_1^2}{R^2}\right] \bigg\}  t_1 G\nonumber\\
&+\bigg\{ \frac{3v^{2/3}}{2 \varepsilon^{1/3}}\left[ (\alpha-1) \frac{|\mathscr G_x|}{v}+ \alpha (p-1)  |\mathscr G_{x v}|
+\alpha (p-1)  \ell_2/2 + \alpha (p-1) \underline k_1 \ell_1\right] ^{4/3} \nonumber\\
&+\beta \mathscr G_v - \alpha (p-1) \Delta_\phi \mathscr G^x
 +\frac{\beta}{\alpha}\left[ \alpha'-\frac{\beta}{b \alpha} \right] - \beta'
+ \alpha^2 (p-1)n[(\underline k_1 + \overline k_1)^2+2k_2] \nonumber\\
&+\frac{b \alpha^2}{ 4 (\alpha -1)^2-2 \varepsilon b \alpha^2} \bigg[ 2(p-1)v[(m-1) k+k_2]
 +2 (\alpha -1) \overline k_1 - \frac{2(\alpha -1)}{b \alpha^2} \beta
 - \frac{\alpha'}{\alpha}  \nonumber\\
&+ (\alpha -1)\left[\frac{\mathscr G}{v} -\mathscr G_v \right]
- \alpha(p-1) v\mathscr G_{vv} \bigg]_+^2 \bigg\} t_1^2 \eta ^2.
\end{align}

Referring to \eqref{eq-2.49-pre}, let us now proceed by denoting the coefficient of $t_1 G$ by ${\bf I}$ and 
the coefficient of $t_1^2 \eta^2=t_1^2 \eta^2 G^0$ by ${\bf II}$, that is, 
\begin{align} \label{EQ-4.23.A}
{\bf I} = &~ \frac{b \alpha^2 p^2 v}{2(\alpha -1)} \frac{c_1 ^2}{R^2}
+\left[\mathscr G_v+\frac{\alpha'}{\alpha} -\frac{2 \beta}{b \alpha^2} +\frac{1}{t_1} \right] \eta+ 2c_1 \underline k_1\nonumber\\
&+(p-1)v \left[\frac{1}{R^2}[c_2 +(m-1) c_1(1+R \sqrt{k})] 
+ 2 \frac{c_1^2}{R^2}\right],
\end{align}
and 
\begin{align} \label{EQ-4.25}
{\bf II} = &~ \frac{3 v^{2/3}}{2\varepsilon^{1/3}} 
\left[ (\alpha-1) \frac{|\mathscr G_x|}{v}+ \alpha (p-1)  |\mathscr G_{x v}|
+\alpha (p-1)  \ell_2/2 + \alpha (p-1) \underline k_1 \ell_1\right]^{4/3}\nonumber\\
&+\beta \mathscr G_v - \alpha (p-1) \Delta_\phi \mathscr G^x 
+ \frac{\beta}{\alpha}\left[ \alpha' -\frac{\beta}{b \alpha} \right]
+ \alpha^2 (p-1)n[(\underline k_1 + \overline k_1)^2+2k_2] \nonumber\\
&- \beta'  +\frac{b \alpha^2}{ 4 (\alpha -1)^2-2 \varepsilon b \alpha^2} 
\bigg[ 2(p-1)v[(m-1) k+k_2] +2 (\alpha -1) \overline k_1\nonumber\\
& - \frac{2(\alpha -1)}{b \alpha^2} \beta
- \frac{\alpha'}{\alpha} + (\alpha -1)\left[\frac{\mathscr G}{v} -\mathscr G_v \right]
- \alpha(p-1) v\mathscr G_{vv} \bigg]_+^2.
\end{align}
Then inequality \eqref{eq-2.49-pre} can be re-written in the following simplified form  
\begin{equation} \label{eq-2.49}
0 \le -\frac{1}{b \alpha^2} G^2 + t_1 {\bf I} G + t_1^2 \eta^2 {\bf II}.
\end{equation}

Denoting by $M = M(2R) = \sup v$ [note that the supremum is taken over $Q_{2R, T}$], and making 
use of \eqref{EQ-eq-K}--\eqref{EQ-eq-E} and \eqref{V-A-1}--\eqref{V-A-4} 
[with $2R$ replacing $R$], it is easily seen from \eqref{EQ-4.23.A} and \eqref{EQ-4.25} 
that we have the bounds
\begin{align}\label{eq2.51}
{\bf I} \le \mathsf A &= \frac{1}{t_1} + \mathsf A^{\star} \\
&= \frac{1}{t_1} +\mu_1(2R)+(p-1)M \left[\frac{1}{R^2}[c_2 +(m-1) c_1(1+R \sqrt{k})] 
+ 2 \frac{c_1^2}{R^2}\right], \nonumber 
\end{align}
and
\begin{align}\label{eq-6.25}
{\bf II} \le \mathsf B = \mu_2^{4/3}(2R)+ \mu_3(2R) +\mu_4^2(2R).
\end{align}
As a result, it follows from \eqref{eq-2.49} that 
\begin{align}
0 \le -\frac{1}{b \alpha^2} G^2 + t_1 \mathsf A G + t_1 ^2 \mathsf B.
\end{align}
We now make use of the implication $-\mathsf a x^2 +\mathsf bx +\mathsf c \ge 0 \implies x \le \mathsf b/\mathsf a +\sqrt{\mathsf c/\mathsf a}$ 
(with ${\mathsf a}>0$, ${\mathsf b}, {\mathsf c} \ge 0$) to deduce that 
\begin{align}\label{eq-6.31}
G(x_1, t_1) \le b t_1 \alpha^2 (t_1) \mathsf A + t_1 \alpha(t_1) \sqrt{b \mathsf B}.
\end{align}
Since $\eta \equiv 1$ for $d(x,x_0, T_1) \le R$ and $(x_1, t_1)$ is the point where $G= t \eta F$ attains its maximum 
on $Q_{2R,T_1}$ we have
\begin{align}
G(x, T_1)= T_1 F(x, T_1) = [t \eta F](x, T_1)  \le [t \eta F] (x_1, t_1) = t_1 [\eta F](x_1, t_1)= G(x_1, t_1).
\end{align}
Hence making note of the bound \eqref{eq-6.31} and dividing through by $T_1 >0$ gives
\begin{align}
F(x,T_1) \le \frac{1}{T_1} G(x_1, t_1) 
&\le b \alpha^2(t_1) \frac{t_1}{T_1} \mathsf A + \alpha(t_1) \frac{t_1}{T_1} \sqrt{b \mathsf B} \nonumber\\
&= b \alpha^2(t_1) \frac{t_1}{T_1} \left(\frac{1}{t_1} + \mathsf A^\star \right) + \alpha(t_1) \frac{t_1}{T_1} \sqrt{b \mathsf B},
\end{align}
and so
\begin{align}
F(x,T_1) &\le \frac{b \alpha^2(t_1)}{T_1} +  \frac{t_1}{T_1} [b \alpha^2(t_1) \mathsf A^\star + \alpha(t_1) \sqrt{b \mathsf B}]\nonumber\\
&\le  \frac{b \alpha^2(T_1)}{T_1} + [b \alpha^2(T_1) \mathsf A^\star + \alpha(T_1) \sqrt{b \mathsf B}]. 
\end{align}
Therefore noting $0 < t_1 \le T_1$, we can write
\begin{align} 
\left[\frac{|\nabla v|^2}{v} -\alpha \frac{\partial_t v}{v} + \alpha \frac{\mathscr G}{v} -\beta \right] (x, T_1)
\le&~ \frac{b \alpha^2(T_1)}{T_1} + b \alpha^2(T_1) \mathsf A^\star + \alpha(T_1) \sqrt{b \mathsf B}.
\end{align}
The arbitrariness of $T_1 >0$ gives the desired conclusion for any $0< t \le T$. $\hfill \square$

\section{Proof of the Harnack inequality in Corollary $\ref{paraHar}$}
\label{sec5}

As this is a global Harnack inequality the idea is to integrate 
the global estimate \eqref{estimate-5.4-global} along suitable space-time curves 
in $\mathscr M \times [0, T]$. To this end, let us more generally consider first the 
case $\alpha=\alpha(t)$ and proceed by re-writing the latter inequality as
\begin{align} \label{newhar}
\frac{|\nabla v|^2}{v} - \alpha(t) \frac{\partial_t v}{v} - \mathsf h (t) \le 0,
\end{align}
or
\begin{align}
- \frac{\partial_t v}{v} \le \frac{1}{\alpha(t)} \left( \mathsf h (t) -\frac{|\nabla v|^2}{v}\right),
\end{align}
where the choice of $\mathsf h(t)$ will be specified below [see \eqref{V-A-H}, \eqref{V-A-h}].
Now upon setting $f= \log v$ and $\underline M = \inf v$ where the infimum is taken over $\mathscr M \times[0, T]$ we have 
\begin{align} \label{Heq8.1}
- \partial_t f = - \frac{\partial_t v}{v} 
&\le \frac{1}{\alpha(t)} \left( \mathsf h (t) -\frac{|\nabla v|^2}{v}\right)
=  \frac{1}{\alpha(t)} \left( \mathsf h (t) - e^f |\nabla f|^2 \right)\nonumber\\
&\le \frac{1}{\alpha(t)} (\mathsf h (t) - \underline M |\nabla f|^2).
\end{align}

Let us fix $x_1, x_2 \in \mathscr M$ and $0<t_1<t_2<T$.
Suppose $\gamma \in \mathscr{C}^1( [t_1,t_2]; \mathscr M)$ is an arbitrary curve with $\gamma(t_1) = x_1$ 
and $\gamma(t_2) = x_2$. Using \eqref{Heq8.1} and writing $\dot \gamma = d\gamma/dt $ it is seen that 
\begin{align}\label{eq-6.39-v}
f(x_1, t_1) - f(x_2, t_2) 
&= \int_{t_2}^{t_1}\frac{d}{dt} [f(\gamma(t),t)] \, dt 
= \int_{t_2}^{t_1} [\langle \nabla f(\gamma(t),t), \dot \gamma (t) \rangle + \partial_t f(\gamma(t),t)] \, dt\nonumber \\
&= \int_{t_1}^{t_2} [-\langle \nabla f(\gamma(t),t), \dot \gamma (t) \rangle - \partial_t f(\gamma(t),t)] \, dt\nonumber \\
&\le \int_{t_1}^{t_2} \left(|\nabla f(\gamma(t),t)|| \dot \gamma (t)| 
+\frac{1}{\alpha(t)} [\mathsf h (t) -\underline M |\nabla f (\gamma(t),t)|^2] \right) \, dt \nonumber\\
& \le \int_{t_1}^{t_2} \left( -\frac{\underline M}{\alpha(t)}  |\nabla f(\gamma(t),t)|^2 + |\nabla f(\gamma(t),t)|| \dot \gamma (t)|\right) \, dt
+  \int_{t_1}^{t_2} \frac{\mathsf h (t)}{\alpha(t)} \, dt \nonumber\\
& \le \int_{t_1}^{t_2} \frac{| \dot \gamma (t)|^2}{4 \underline M} \alpha (t) \, dt
+  \int_{t_1}^{t_2} \frac{\mathsf h (t)}{\alpha(t)} \, dt,
\end{align}
where we have used $-\mathsf a x^2 + \mathsf b x \le \mathsf b^2 /(4 \mathsf a)$. [Note that in the event $\underline{M}=0$ 
by interpretting $1/\underline{M}$ as infinity the inequalities \eqref{eq-6.39-v} and \eqref{eq-6.43-a} below are trivially true.]
Note also that in \eqref{eq-6.39-v} all inner products are with respect to the metric $g(t)$ where $t_1\le t \le t_2$.
Exponentiating the above inequality gives 
\begin{align}\label{eq-6.43-a}
v(x_1, t_1) \le v(x_2, t_2)~ {\rm exp} \left[ \frac{1}{4 \underline M} \int_{t_1}^{t_2} | \dot \gamma (t)|^2 \alpha (t) \, dt
+  \int_{t_1}^{t_2} \frac{\mathsf h (t)}{\alpha(t)} \, dt\right].
\end{align}
Now to complete the proof of the Harnack inequality \eqref{Eq5.62Har}, 
we first write 
\begin{align}\label{V-A-h}
\mathsf h(t) = \frac{b \alpha^2}{t} + \mathsf H,
\end{align}
where $\mathsf H$ is as in \eqref{V-A-H}.
Then in view of the global estimate \eqref{estimate-5.4-global} we have the inequality \eqref{newhar} 
with the choice of $\mathsf h(t)$ as in \eqref{V-A-h}. 
The conclusion of the corollary now follows by reparameterising $\gamma$, 
changing variables in the integral \eqref{eq-6.43-a} and noting that $\alpha>1$ is constant. 
For the sake of clarity note that here as $|\dot \gamma|^2 = |\dot \gamma|^2_{g}$ 
the stated change of variables will also apply to the arguments of $g$. \hfill $\square$

\section{Proof of the local estimate in Theorem \ref{Thrm-2.5}}  
\label{sec6}

Starting with the description of the evolution $\mathscr L^p_v [F] = [\partial_t -(p-1) v \Delta_\phi] F$ in \eqref {eq-5.27}, 
using the curvature lower bound ${\mathscr Ric}_\phi^m (g) \ge -(m-1)k g$ and the inequality 
\begin{align}
\frac{(\Delta_\phi  v)^2}{m} - \alpha^2 (t) |h|^2 
& \le |\nabla \nabla  v|^2 + \frac{\langle \nabla \phi , \nabla v \rangle^2}{m-n} - \alpha^2 (t) |h|^2 \nonumber \\
& \le 2|\nabla \nabla v|^2 + 2\frac{\langle \nabla \phi , \nabla v \rangle^2}{m-n}
- 2\alpha(t) \langle h , \nabla \nabla v \rangle,  
\end{align}
we can write 
\begin{align}
\mathscr L^p_v [F] \le& - \frac{1+m(p-1)}{m(p-1)}[(p-1)\Delta_\phi v]^2
+ 2p \langle \nabla F, \nabla v \rangle
+2(p-1)v (m-1) k \frac{|\nabla v|^2}{v} \nonumber\\
&+ [1-\alpha (t)] \left[ \frac{\partial_t v}{v} - \frac{\mathscr G(t,x,v)}{v} \right]^2
+ (p-1) \alpha^2(t) |h|^2 + \frac{2}{v}[\alpha(t) -1] h ( \nabla v, \nabla v)\nonumber\\
&+ \alpha(t) (p-1) \langle 2 {\rm div} h - \nabla ({\rm Tr}_g h),\nabla v \rangle
+ \alpha(t) (p-1) \langle \nabla \partial_t \phi , \nabla v \rangle-\beta' (t)\nonumber\\
&- 2 \alpha(t) (p-1)h (\nabla \phi , \nabla v) 
+\frac{2}{v}[1-\alpha(t)] \langle \nabla \mathscr G(t,x,v), \nabla v \rangle
+ \alpha'(t) \frac{\mathscr G(t,x,v)}{v}\nonumber\\
& -\alpha'(t) \frac{\partial_t v}{v} - \alpha (t)(p-1) \Delta_\phi \mathscr G(t,x,v) 
+[\alpha(t)-1] \frac{|\nabla v|^2}{v} \frac{\mathscr G(t,x,v)}{v}.
\end{align}
In particular, by applying $\mathscr L^p_v$ to the quotient $F/\alpha$ it follows from the above that  
\begin{align}
\mathscr L^p_v [F/\alpha] =&~ -\frac{\alpha'}{\alpha^2} F + \frac{1}{\alpha} \mathscr L^p_v [F]\nonumber\\
\le& -\frac{\alpha'}{\alpha^2} \frac{|\nabla v|^2}{v} + \frac{\alpha'}{\alpha}\frac{\partial_t v}{v} 
-\frac{\alpha'}{\alpha} \frac{\mathscr G}{v} +\frac{\alpha'}{\alpha^2}\beta 
- \frac{1}{b\alpha}[(p-1)\Delta_\phi v]^2+ \frac{2p}{\alpha} \langle \nabla F, \nabla v \rangle\nonumber\\
&+\frac{2(p-1)}{\alpha}v (m-1) k \frac{|\nabla v|^2}{v} + (p-1) \alpha |h|^2
+ \frac{1-\alpha}{\alpha}\left[ \frac{\partial_t v}{v} - \frac{\mathscr G}{v} \right]^2\nonumber\\
&-\frac{\alpha'}{\alpha} \frac{\partial_t v}{v} +\frac{2}{v}\frac{\alpha -1}{\alpha} h ( \nabla v, \nabla v)
+(p-1) \langle 2 {\rm div} h - \nabla ({\rm Tr}_g h),\nabla v \rangle\nonumber\\
&+ (p-1) \langle \nabla \partial_t \phi , \nabla v \rangle-\frac{\beta'}{\alpha}
- 2 (p-1)h (\nabla \phi , \nabla v) +\frac{2}{v}\frac{1-\alpha}{\alpha} \langle \nabla \mathscr G, \nabla v \rangle\nonumber\\
&+ \frac{\alpha'}{\alpha} \frac{\mathscr G}{v} - (p-1) \Delta_\phi \mathscr G 
+\frac{\alpha-1}{\alpha} \frac{|\nabla v|^2}{v} \frac{\mathscr G}{v}.
\end{align}

Using \eqref{eq-5.56}, \eqref{eq-5.57}, \eqref{sarkhoob-pb1}, the bound 
$\alpha |h|^2  \le n \alpha (\underline k_1 + \overline k_1)^2$ along with \eqref{EQ-eq-3.45}, \eqref{EQ-eq-3.46} in the form of 
$\langle  2 {\rm div} h - \nabla ({\rm Tr}_g h),\nabla v \rangle \le 2 n k_2 + 2 k_2 |\nabla v|^2$ 
and \eqref{EQ-eq-3.47}, \eqref{EQ-eq-3.48} results in 
\begin{align}
\mathscr L^p_v [F/\alpha] 
\le& -\frac{F^2}{b\alpha^3} + \frac{2(1-\alpha)}{b\alpha^2}
\frac{|\nabla v|^2}{v} \frac{F}{\alpha}
+ \left[\mathscr G_v-\frac{2\beta}{b\alpha^2}\right]\frac{F}{\alpha}
+ 2p\langle \nabla (F/\alpha), \nabla v \rangle \\
&-\frac{(1-\alpha)^2}{b\alpha^3}\frac{|\nabla v|^4}{v^2}
+\bigg[2(p-1)v \left(\frac{(m-1) k}{\alpha} +k_2\right)
+ \frac{2(\alpha -1)}{\alpha} \overline k_1\nonumber\\
&+\frac{\alpha-1}{\alpha} \left(\frac{\mathscr G}{v}-  \mathscr G_v\right)
-(p-1)v \mathscr G_{vv}+\frac{2\beta(1-\alpha)}{b\alpha^3}
-\frac{\alpha'}{\alpha^2}\bigg] \frac{|\nabla v|^2}{v}\nonumber\\
&+\left[\frac{2 (\alpha-1)}{\alpha} \frac{|\mathscr G_x|}{v}+2 (p-1) |\mathscr G_{x v}|
+(p-1) \ell_2 + 2 (p-1) \underline k_1 \ell_1 \right] |\nabla v| \nonumber\\
&+ (p-1) n[\alpha (\underline k_1 + \overline k_1)^2 + 2 k_2]
- (p-1) \Delta_ \phi \mathscr G^x + \frac{\beta}{\alpha} \mathscr G_v
-\frac{\beta}{\alpha^2}\left(\frac{\beta}{b\alpha}-\alpha' \right)
-\frac{\beta'}{\alpha}.\nonumber
\end{align}

Let $G = t \eta F/\alpha$ where $\eta$ is the same cut-off function as before. 
Then with $(x_1, t_1)$ the maximum point 
for $G$ over $Q_{2R, T_1} \subset Q_{2R, T}$ we have
at this maximum point $\nabla (F/\alpha)= - F\nabla \eta/(\alpha \eta)$
and so $\langle \nabla v, \nabla (F/\alpha) \rangle \le |\nabla v| |\nabla (F/\alpha)| 
= |\nabla v| (|\nabla \eta|/\eta) F/\alpha$. Thus  
\begin{align}
0 \le&~ \mathscr L^p_v [G] = t_1 \eta  \mathscr L^p_v[F/\alpha] 
+2  t_1 (p-1) v\frac{|\nabla \eta |^2}{\eta} \frac{F}{\alpha} 
+ t_1 \frac{F}{\alpha} \mathscr L^p_v [\eta] + \eta \frac{F}{\alpha}\\
\le &~2  t_1 (p-1) v\frac{|\nabla \eta |^2}{\eta} \frac{F}{\alpha} 
+ t_1 \frac{F}{\alpha} \mathscr L^p_v [\eta] + \eta \frac{F}{\alpha} 
-t_1 \eta \frac{F^2}{b\alpha^3} + \frac{2(1-\alpha)}{b\alpha^2}
\frac{|\nabla v|^2}{v} \frac{ t_1 \eta F}{\alpha}\nonumber\\
&+ \left[\mathscr G_v-\frac{2\beta}{b\alpha^2}\right]t_1 \eta \frac{F}{\alpha}
+ 2p t_1 \eta |\nabla v| \frac{|\nabla \eta|}{\eta} \frac{F}{\alpha} 
-\frac{(1-\alpha)^2}{b\alpha^3} t_1 \eta \frac{|\nabla v|^4}{v^2}\nonumber\\
&+\bigg[2(p-1)v \left(\frac{(m-1) k}{\alpha} +k_2\right)
+ \frac{2(\alpha -1)}{\alpha} \overline k_1
+\frac{\alpha-1}{\alpha} \left(\frac{\mathscr G}{v}-  \mathscr G_v\right)\nonumber\\
&-(p-1)v \mathscr G_{vv}+\frac{2\beta(1-\alpha)}{b\alpha^3}
-\frac{\alpha'}{\alpha^2}\bigg] t_1 \eta \frac{|\nabla v|^2}{v}
+ t_1 \eta |\nabla v| \times \nonumber\\
&\times \left[\frac{2 (\alpha-1)}{\alpha} \frac{|\mathscr G_x|}{v}+2 (p-1) |\mathscr G_{x v}|
+(p-1) \ell_2 + 2 (p-1) \underline k_1 \ell_1 \right] \nonumber\\
&+ \bigg[(p-1) n[\alpha (\underline k_1 + \overline k_1)^2 + 2 k_2]
- (p-1) \Delta_ \phi \mathscr G^x + \frac{\beta}{\alpha} \mathscr G_v
-\frac{\beta}{\alpha^2}\left(\frac{\beta}{b\alpha}-\alpha' \right)
-\frac{\beta'}{\alpha}\bigg] t_1 \eta .\nonumber
\end{align}
Using the identity $G = t_1 \eta F/\alpha$ and multiplying through by $\eta$ then gives
\begin{align}\label{EQ-3.46}
0 \le &-\frac{1}{b \alpha}\frac{G^2}{t_1}
+ \frac{2(1-\alpha)}{b\alpha^2}\frac{|\nabla v|^2}{v} \eta G
+ 2p\sqrt v \frac{|\nabla v|}{\sqrt v} \frac{|\nabla \eta|}{\eta} \eta G
-\frac{(1-\alpha)^2}{b\alpha^3} t_1 \eta^2 \frac{|\nabla v|^4}{v^2}\\
&+ \left[\frac{\partial_t \eta}{\eta}-(p-1)v \frac{\Delta_\phi \eta}{\eta}
+2(p-1) v\frac{|\nabla \eta |^2}{\eta^2}+ \frac{1}{t_1}
-\frac{2\beta}{b\alpha^2} +\mathscr G_v\right] \eta G \nonumber\\
&+\bigg[2(p-1)v \left(\frac{(m-1) k}{\alpha} +k_2\right)
+ \frac{2(\alpha -1)}{\alpha} \overline k_1\nonumber\\
&+\frac{\alpha-1}{\alpha} \left(\frac{\mathscr G}{v}-  \mathscr G_v\right)
-(p-1)v \mathscr G_{vv}+\frac{2\beta(1-\alpha)}{b\alpha^3}
-\frac{\alpha'}{\alpha^2}\bigg] t_1 \eta^2 \frac{|\nabla v|^2}{v}\nonumber\\
&+\left[\frac{2 (\alpha-1)}{\alpha} \frac{|\mathscr G_x|}{v}+2 (p-1) |\mathscr G_{x v}|
+(p-1) \ell_2 + 2 (p-1) \underline k_1 \ell_1 \right] t_1 \eta^2 |\nabla v| \nonumber\\
&+ \left[(p-1) n[\alpha (\underline k_1 + \overline k_1)^2 + 2 k_2]
- (p-1) \Delta_ \phi \mathscr G^x + \frac{\beta}{\alpha} \mathscr G_v
-\frac{\beta}{\alpha^2}\left(\frac{\beta}{b\alpha}-\alpha' \right)
-\frac{\beta'}{\alpha}\right] t_1 \eta^2.\nonumber
\end{align}
Using \eqref{Lpv-Ali-Added}, \eqref{Eq5.35}, \eqref{eq5.68} and Young's inequality 
as in \eqref{eq5.72} but now in the form 
\begin{align}
\bigg\{ &\frac{(\alpha-1)}{\alpha} \frac{|\mathscr G_x|}{v}+ (p-1)  |\mathscr G_{x v}|+ 
(p-1)  \ell_2/2+(p-1) \underline k_1 \ell_1 \bigg\}|\nabla v| \\
&\le \frac{3v^{2/3}}{4\varepsilon^{1/3}} \left\{ \frac{(\alpha-1)}{\alpha} \frac{|\mathscr G_x|}{v}
+(p-1)  |\mathscr G_{x v}| +(p-1)  \ell_2/2 
+ (p-1) \underline k_1 \ell_1 \right\}^{4/3} + \varepsilon\frac{|\nabla v|^4}{4v^2}, \nonumber
\end{align}
we can write 
\begin{align} \label{EQ-eq-6.9}
0 \le &-\frac{1}{b \alpha}\frac{G^2}{t_1}
+\frac{b \alpha^2 p^2 v}{2(\alpha-1)} \frac{c_1 ^2}{R^2} G
-\left[\frac{(1-\alpha)^2}{b\alpha^3} -\frac{\varepsilon}{2}\right]t_1 \eta^2 \frac{|\nabla v|^4}{v^2}\\
&+\left\{ (p-1)v \left[ \frac{1}{R^2}[c_2 +(m-1) c_1(1+R \sqrt{k})]
+ 2 \frac{c_1^2}{R^2}\right]+ 2c_1 \underline k_1 \right\} G\nonumber\\
&+ \left[\frac{1}{t_1}-\frac{2\beta}{b\alpha^2} +\mathscr G_v\right] \eta G 
+\bigg[2(p-1)v \left(\frac{(m-1) k}{\alpha} +k_2\right)
+ \frac{2(\alpha -1)}{\alpha} \overline k_1\nonumber\\
&+\frac{\alpha-1}{\alpha} \left(\frac{\mathscr G}{v}-  \mathscr G_v\right)
-(p-1)v \mathscr G_{vv}+\frac{2\beta(1-\alpha)}{b\alpha^3}
-\frac{\alpha'}{\alpha^2}\bigg] t_1 \eta^2 \frac{|\nabla v|^2}{v}\nonumber\\
&+\frac{3v^{2/3}}{2\varepsilon^{1/3}} \left\{ \frac{(\alpha-1)}{\alpha} \frac{|\mathscr G_x|}{v}
+(p-1)  |\mathscr G_{x v}| +(p-1)  \ell_2/2 
+ (p-1) \underline k_1 \ell_1 \right\}^{4/3} t_1 \eta^2\nonumber\\
&+ \left[(p-1) n[\alpha (\underline k_1 + \overline k_1)^2 + 2 k_2]
- (p-1) \Delta_ \phi \mathscr G^x + \frac{\beta}{\alpha} \mathscr G_v
-\frac{\beta}{\alpha^2}\left(\frac{\beta}{b\alpha}-\alpha' \right)
-\frac{\beta'}{\alpha}\right] t_1 \eta^2.\nonumber
\end{align}

Again, by using the inequality $-\mathsf ax^2 +\mathsf bx \le \mathsf b^2 /(4\mathsf a)$ 
with $x$ being $|\nabla v|^2/v$ and $\mathsf a > 0$, and noting $0<\varepsilon<2(1-\alpha)^2/(b \alpha^3)$, 
it is evident that
\begin{align}
-\frac{2(1-\alpha)^2- \varepsilon b \alpha^3}{2b\alpha^3}& t_1 \eta^2 \frac{|\nabla v|^4}{v^2}
+\bigg[2(p-1)v \left(\frac{(m-1) k}{\alpha} +k_2\right)
+ \frac{2(\alpha -1)}{\alpha} \overline k_1\nonumber\\
&+\frac{\alpha-1}{\alpha} \left(\frac{\mathscr G}{v}-  \mathscr G_v\right)
-(p-1)v \mathscr G_{vv}+\frac{2\beta(1-\alpha)}{b\alpha^3}
-\frac{\alpha'}{\alpha^2}\bigg] t_1 \eta^2 \frac{|\nabla v|^2}{v}\nonumber\\
&\le \frac{b \alpha^3 t_1 \eta^2}{ 4 (1-\alpha)^2-2\varepsilon b \alpha^3} 
\bigg[2(p-1)v \left(\frac{(m-1) k}{\alpha} +k_2\right)
+ 2\frac{\alpha -1}{\alpha} \overline k_1\nonumber\\
&+\frac{\alpha-1}{\alpha} \left(\frac{\mathscr G}{v}-  \mathscr G_v\right)
-(p-1)v \mathscr G_{vv}+\frac{2\beta(1-\alpha)}{b\alpha^3}
-\frac{\alpha'}{\alpha^2}\bigg]_+^2.
\end{align}
Substituting the above in \eqref{EQ-eq-6.9} and multiplying through by $t_1>0$ gives
\begin{align} \label{EQ-eq-6.11}
0 \le &-\frac{1}{b \alpha}G^2
+\bigg\{\frac{b \alpha^2 p^2 v}{2(\alpha-1)} \frac{c_1 ^2}{R^2}
+(p-1)v \left[ \frac{1}{R^2}[c_2 +(m-1) c_1(1+R \sqrt{k})]
+ 2 \frac{c_1^2}{R^2}\right]+ 2c_1 \underline k_1 \nonumber\\
&+ \left[\frac{1}{t_1}-\frac{2\beta}{b\alpha^2} + \mathscr G_v\right] \eta \bigg\} t_1 G 
+\bigg\{\frac{b \alpha^3}{ 4 (1-\alpha)^2-2\varepsilon b \alpha^3} 
\bigg[2(p-1)v \left(\frac{(m-1) k}{\alpha} +k_2\right)\nonumber\\
&+\frac{2(\alpha -1)}{\alpha} \overline k_1+\frac{\alpha-1}{\alpha}
 \left(\frac{\mathscr G}{v}-  \mathscr G_v\right)
-(p-1)v \mathscr G_{vv}+\frac{2\beta(1-\alpha)}{b\alpha^3}
-\frac{\alpha'}{\alpha^2}\bigg]_+^2 \\
&+\frac{3v^{2/3}}{2\varepsilon^{1/3}} \left[ \frac{(\alpha-1)}{\alpha} \frac{|\mathscr G_x|}{v}
+(p-1)  |\mathscr G_{x v}| +(p-1)  \ell_2/2 
+ (p-1) \underline k_1 \ell_1 \right]^{4/3}\nonumber\\
&+ \left[(p-1) n[\alpha (\underline k_1 + \overline k_1)^2 + 2 k_2]
- (p-1) \Delta_ \phi \mathscr G^x + \frac{\beta}{\alpha} \mathscr G_v
-\frac{\beta}{\alpha^2}\left(\frac{\beta}{b\alpha}-\alpha' \right)
-\frac{\beta'}{\alpha}\right]\bigg\} t_1^2 \eta^2.\nonumber
\end{align}

Referring to \eqref{EQ-eq-6.11}, let us now proceed by denoting the coefficient of $t_1 G$ by ${\bf I}$ and 
the coefficient of $t_1^2 \eta^2=t_1^2 \eta^2 G^0$ by ${\bf II}$, that is, 

\begin{align} \label{EQ-eq-6.12.A}
{\bf I} = &~ \frac{b \alpha^2 p^2 v}{2(\alpha -1)} \frac{c_1 ^2}{R^2}
+\left[ \mathscr G_v-\frac{2 \beta}{b \alpha^2} +\frac{1}{t_1} \right] \eta+ 2c_1 \underline k_1\nonumber\\
&+(p-1)v \left[\frac{1}{R^2}[c_2 +(m-1) c_1(1+R \sqrt{k})] 
+ 2 \frac{c_1^2}{R^2}\right],
\end{align}
and 
\begin{align} \label{EQ-eq-6.13.A}
{\bf II} = &~ \frac{3v^{2/3}}{2\varepsilon^{1/3}} 
\left[ \frac{\alpha-1}{\alpha} \frac{|\mathscr G_x|}{v}+ (p-1)  |\mathscr G_{x v}|
+ (p-1)  \ell_2/2 + (p-1) \underline k_1 \ell_1\right]^{4/3}\nonumber\\
&+\frac{\beta}{\alpha} \mathscr G_v - (p-1) \Delta_\phi \mathscr G^x 
+ \frac{\beta}{\alpha^2}\left[ \alpha' -\frac{\beta}{b \alpha} \right] - \frac{\beta'}{\alpha} 
+ (p-1)n[\alpha (\underline k_1 + \overline k_1)^2+2k_2] \nonumber\\
&+\frac{b \alpha^3}{ 4 (1-\alpha)^2-2\varepsilon b \alpha^3} 
\bigg[2(p-1)v \left(\frac{(m-1) k}{\alpha} +k_2\right)+ \frac{2(\alpha -1)}{\alpha} \overline k_1\nonumber\\
&+\frac{\alpha-1}{\alpha} \left(\frac{\mathscr G}{v}-  \mathscr G_v\right)
-(p-1)v \mathscr G_{vv}+\frac{2\beta(1-\alpha)}{b\alpha^3}
-\frac{\alpha'}{\alpha^2}\bigg]_+^2.
\end{align}
Then inequality \eqref{EQ-eq-6.11} can be re-written in the following simplified form  
\begin{equation} \label{eq-6.13}
0 \le -\frac{1}{b \alpha} G^2 + t_1 {\bf I} G + t_1^2 \eta^2 {\bf II}.
\end{equation}
Denoting by $M = M(2R) = \sup v$ and making 
use of \eqref{EQ-eq-K}, \eqref{EQ-eq-L}, \eqref{EQ-eq-E}, \eqref{EQ-eq-M.n}--\eqref{EQ-eq-E.n} and \eqref{V-A-1.n}--\eqref{V-A-4.n} 
[with $2R$ replacing $R$], it is easily seen from \eqref{EQ-eq-6.12.A} and \eqref{EQ-eq-6.13.A} 
that we have the bounds
\begin{align}\label{eq2.51}
{\bf I} \le \mathsf A &= \frac{1}{t_1} + \mathsf A^{\star} \\
&= \frac{1}{t_1} +\lambda_1(2R)+(p-1)M \left[\frac{1}{R^2}[c_2 +(m-1) c_1(1+R \sqrt{k})] 
+ 2 \frac{c_1^2}{R^2}\right], \nonumber 
\end{align}
and
\begin{align}\label{eq-6.25}
{\bf II} \le \mathsf B = \lambda_2^{4/3}(2R)+ \lambda_3(2R) +\lambda_4^2(2R).
\end{align}
As a result, it follows from \eqref{eq-6.13} that 
\begin{align}
0 \le -\frac{1}{b \alpha} G^2 + t_1 \mathsf A G + t_1 ^2 \mathsf B.
\end{align}
We now make use of the implication $-\mathsf a x^2 +\mathsf bx +\mathsf c \ge 0 \implies x \le \mathsf b/\mathsf a +\sqrt{\mathsf c/\mathsf a}$ 
(with ${\mathsf a}>0$, ${\mathsf b}, {\mathsf c} \ge 0$) to deduce that 
\begin{align} \label{New-eq-6.31}
G(x_1, t_1) \le b \alpha(t_1) t_1 \mathsf A + t_1 \sqrt{b \alpha(t_1) \mathsf B}.
\end{align}
Since $\eta \equiv 1$ for $d(x,x_0, T_1) \le R$ and $(x_1, t_1)$ is the point where $G= t \eta F/\alpha$ attains its maximum 
on $Q_{2R,T_1}$ we have
\begin{align}
G(x, T_1)= \frac{T_1 F(x, T_1)}{\alpha(T_1)} 
= \frac{[t \eta F](x, T_1)}{\alpha(T_1)} 
\le \frac{[t \eta F](x_1, t_1)}{\alpha(t_1)} = G(x_1, t_1).
\end{align}
Hence making note of the bound \eqref{New-eq-6.31} and dividing through by $T_1 >0$ gives
\begin{align}
\frac{F(x,T_1)}{\alpha(T_1)} \le \frac{1}{T_1} G(x_1, t_1) 
&\le b \alpha(t_1) \frac{t_1}{T_1} \mathsf A 
+ \frac{t_1}{T_1} \sqrt{b \alpha(t_1) \mathsf B} \nonumber \\
&= b \alpha(t_1) \frac{t_1}{T_1} \left(\frac{1}{t_1} 
+ \mathsf A^\star \right) + \frac{t_1}{T_1} \sqrt{b \alpha(t_1) \mathsf B},
\end{align}
and so in view of $0< t_1 \le T_1$ and monotonicity of $\alpha$ this gives
\begin{align}
\frac{F(x,T_1)}{\alpha(T_1)} &\le \frac{b \alpha(t_1)}{T_1} +  \frac{t_1}{T_1} 
[b \alpha(t_1) \mathsf A^\star + \sqrt{b \alpha(t_1) \mathsf B}] \nonumber \\
&\le \frac{b \alpha(T_1)}{T_1} + 
[b \alpha(T_1) \mathsf A^\star + \sqrt{b \alpha(T_1) \mathsf B}].
\end{align}
The arbitrariness of $T_1>0$ now gives the desired conclusion.

\section{Proof of the Harnack inequality in Corollary $\ref{corollary-Harnack-global.n}$}
\label{sec7}

Referring to the proof of Corollary \ref{paraHar} in Section \ref{sec5} we start again with the inequality \eqref{eq-6.43-a} and writing 
\begin{align} \label{eq-6.43-a.n}
v(x_1, t_1) \le v(x_2, t_2)~ {\rm exp} \left[ \frac{1}{4 \underline M} \int_{t_1}^{t_2} | \dot \gamma (t)|^2 \alpha (t) \, dt
+  \int_{t_1}^{t_2} \frac{\mathsf h (t)}{\alpha(t)} \, dt\right].
\end{align}
Now to complete the proof of the Harnack inequality \eqref{Eq5.62Har.n}, 
we first write 
\begin{align}\label{V-A-h.n}
\mathsf h(t) = \frac{b \alpha^2}{t} + \mathsf H,
\end{align}
where $\mathsf H$ is as in \eqref{V-A-H.n}.
Then in view of the global estimate \eqref{estimate-2-global} we have the inequality \eqref{newhar} 
with the choice of $\mathsf h(t)$ as in \eqref{V-A-h.n}. 
The conclusion of the corollary now follows by reparameterising $\gamma$, 
changing variables in the integral \eqref{eq-6.43-a.n} and noting that $\alpha>1$ is constant.
\hfill $\square$

\qquad \\
{\bf Acknowledgement.} The authors gratefully acknowledge support from the Engineering and Physical 
Sciences Research Council (EPSRC) through the grant EP/V027115/1. All data is provided in full in the 
results section. Additional data is in the public domain at locations cited in the reference section.

\end{document}